\documentstyle[amscd]{amsart}
\numberwithin{equation}{section}
\theoremstyle{plain}
\newtheorem{thm}{Theorem}[section]
\newtheorem{cor}[thm]{Corollary}
\newtheorem{lem}[thm]{Lemma}
\newtheorem{prop}[thm]{Proposition}

\begin{document}
\title{$C^*$-algebras associated with lambda-synchronizing subshifts and \\
flow equivalence}
\author{Kengo Matsumoto}
\address{ Department of Mathematics, 
Joetsu University of Education,
Joetsu 943-8512,
 Japan}
\email{kengo{@@}juen.ac.jp}
\begin{abstract}
A certain synchronizing property for subshifts called  $\lambda$-synchronization  
yields  $\lambda$-graph systems called the $\lambda$-synchronizing $\lambda$-graph systems
for the subshifts.
The $\lambda$-synchronizing $\lambda$-graph system is a left Fischer cover analogue
for a $\lambda$-synchronizing subshift.
We will study algebraic structure of the $C^*$-algebra
associated with the 
$\lambda$-synchronizing $\lambda$-graph system
and
prove that the stable isomorphism class of the $C^*$-algebra with its Cartan subalgebra 
is invariant under flow equivalence of 
$\lambda$-synchronizing subshifts.
\end{abstract}

\maketitle



\def\Zp{{ {\Bbb Z}_+ }}

\def\U{{ {\cal U} }}
\def\S{{ {\cal S} }}
\def\M{{ {\cal M} }}
\def\P{{ {\cal P} }}
\def\Q{{ {\cal Q} }}
\def\G{{ {\cal G} }}

\def\LLL{{ {\frak L}^{\lambda(\Lambda)} }}
\def\LLTL{{ {\frak L}^{\lambda(\widetilde{\Lambda})} }}

\def\OFL{{ {\cal O}_{\frak L}}}
\def\OLL{{ {\cal O}_{\lambda(\Lambda)}  }}
\def\OLTL{{ {\cal O}_{\lambda(\widetilde{\Lambda})}  }}
\def\ALL{{ {\cal A}_{\lambda(\Lambda)}  }}
\def\ALTL{{ {\cal A}_{\lambda(\widetilde{\Lambda})}  }}
\def\DLL{{ {\cal D}_{\lambda(\Lambda)}  }}
\def\DLTL{{ {\cal D}_{\lambda(\widetilde{\Lambda})}  }}
\def\FL{{{\cal F}_{\frak L}}}
\def\FKL{{ {\cal F}_k^{l} }}
\def\A{{ {\cal A} }}
\def\Ext{{{\operatorname{Ext}}}}
\def\Im{{{\operatorname{Im}}}}
\def\Hom{{{\operatorname{Hom}}}}
\def\Ker{{{\operatorname{Ker}}}}
\def\dim{{{\operatorname{dim}}}}
\def\id{{{\operatorname{id}}}}
\def\OLF{{{\cal O}_{{\frak L}^{Ch(D_F)}}}}
\def\OLN{{{\cal O}_{{\frak L}^{Ch(D_N)}}}}
\def\OLA{{{\cal O}_{{\frak L}^{Ch(D_A)}}}}
\def\LCHDA{{{{\frak L}^{Ch(D_A)}}}}
\def\LCHDF{{{{\frak L}^{Ch(D_F)}}}}
\def\LCHLA{{{{\frak L}^{Ch(\Lambda_A)}}}}
\def\LWA{{{{\frak L}^{W(\Lambda_A)}}}}


\section{Introduction}

Let 
$\Sigma$
be a finite set
with its discrete topology. 
We call it an alphabet and each member of it a symbol or a label.
Let $\Sigma^{\Bbb Z}$, $\Sigma^{\Bbb N}$ 
be the infinite product spaces 
$\prod_{i=-{\infty}}^{\infty}\Sigma_{i}$, 
$\prod_{i=1}^{\infty}\Sigma_i$ 
where 
$\Sigma_{i} = \Sigma$,
 endowed with the product topology respectively.
 The transformation $\sigma$ on $\Sigma^{\Bbb Z}$ 
given by 
$\sigma((x_i)_{i \in {\Bbb Z}}) = (x_{i+1})_{i\in {\Bbb Z}}$ 
for $(x_i)_{i \in {\Bbb Z}} \in \Sigma^{\Bbb Z}$
is called the full shift.
 Let $\Lambda$ be a shift invariant closed subset of $\Sigma^{\Bbb Z}$ i.e. $\sigma(\Lambda) = \Lambda$. 
 The topological dynamical system
   $(\Lambda, \sigma\vert_{\Lambda})$
   is called a subshift or a symbolic dynamical system,
   and simply written as $\Lambda$.
Theory of symbolic dynamical systems 
forms a basic  ingredient
 in the theory of topological dynamical systems
(see \cite{Kit}, \cite{LM} ).

The author has introduced a notion of $\lambda$-graph system,
that is a labeled Bratteli diagram with an additional structure called 
$\iota$-map (\cite{1999DocMath}).
 A $\lambda$-graph system ${\frak L}$ presents a subshift and 
 yields a $C^*$-algebra ${\cal O}_{\frak L}$ (\cite{2002DocMath}).
For a subshift $\Lambda$, 
one may construct a $\lambda$-graph system
${\frak L}^\Lambda$ called the canonical $\lambda$-graph system for $\Lambda$
in a canonical way.
It is a left Krieger cover version for a subshift.
The $C^*$-algebra 
${\cal O}_{{\frak L}^\Lambda}$
for ${\frak L}^\Lambda$
coincides with  the $C^*$-algebra ${\cal O}_\Lambda$
associated with subshift $\Lambda$ (\cite{1997IJM},cf. \cite{CaMa}).
It has been proved that the stable isomorphism class of the $C^*$-algebra
 ${\cal O}_\Lambda$ is invariant under not only topological conjugacy of  $\Lambda$ but also flow equivalence of $\Lambda$,
 so that the K-groups $K_i({\cal O}_\Lambda), i=0,1$
 and the Ext-groups $\Ext^i({\cal O}_\Lambda), i=0,1$
 are invariant under flow equivalence of subshifts (\cite{CK}, \cite{2001K}, \cite{2001ETDS}).
 The latter groups  $\Ext^i({\cal O}_\Lambda), i=0,1$
  have been defined as the Bowen-Franks groups for $\Lambda$ (\cite{2001K}, \cite{2001ETDS}).
For an irreducible sofic shift,
there is another important cover called the (left or right)  Fischer cover.
The (left) Fischer cover is an irreducible labeled graph 
that is  minimal (left)-resolving presentation,
whereas the (left) Krieger cover is not necessarily irreducible.

 In \cite{KM2010},
 a certain synchronizing property for subshifts  
called  $\lambda$-synchronization has been introduced.
The  $\lambda$-synchronizing property  
is weaker than the  usual synchronizing property,
so that  irreducible sofic shifts are $\lambda$-synchronizing 
as well as Dyck sihifts, $\beta$-shifts, Morse shifts, etc.
are $\lambda$-synchronizing.
Many irreducible subshifts have this property.
For a $\lambda$-synchronizing subshift $\Lambda$
there exists a $\lambda$-graph system
called the $\lambda$-synchronizing $\lambda$-graph system
$\LLL$.
The $\lambda$-synchronizing $\lambda$-graph system 
for an irreducible  sofic shift is 
the $\lambda$-graph system associated to the left Fisher cover.
Hence
the $\lambda$-synchronizing $\lambda$-graph systems 
are regarded as the left Fisher cover analogue
for $\lambda$-synchronizing subshifts.

In \cite{2011Pre}, it has been proved that 
the K-groups and the Bowen-Franks groups for 
a  $\lambda$-synchronizing $\lambda$-graph system are invariant under not only topological conjugacy
but also flow equivalence,
so that they yield flow equivalence invariants of $\lambda$-synchronizing subshifts.

In this paper,
we will first study algebraic structure of the $C^*$-algebra
${\cal O}_{\lambda(\Lambda)}$
associated with the $\lambda$-synchronizing $\lambda$-graph system
$\LLL$ for $\Lambda$,
and show that 
if $\Lambda$ is $\lambda$-synchronizingly transitive,
the algebra 
${\cal O}_{\lambda(\Lambda)}$
is simple (Theorem 3.7).
We will next prove that 
the stable isomorphism class of the $C^*$-algebra ${\cal O}_{\lambda(\Lambda)}$
with its Cartan subalgebra
is invariant under flow equivalence of 
$\lambda$-synchronizing subshifts (Theorem 4.17).
As a consequence we have a $C^*$-algebraic proof for the above mentioned fact that 
the K-groups and the Bowen-Franks groups for the  $\lambda$-synchronizing $\lambda$-graph system are invariant under  flow equivalence
(Corollary 4.18).

\section{$\lambda$-synchronizing $\lambda$-graph systems}
Let $\Lambda$ be a subshift over $\Sigma$.
We denote by
$X_{\Lambda} (\subset \Sigma^{\Bbb N})$
the set of all right one-sided sequences appearing in $\Lambda$.
For a natural number $l \in {\Bbb N}$,
we denote by $B_l(\Lambda)$ 
the set of all words appearing in $\Lambda$ 
with length equal to
$l$.
Put
$B_*(\Lambda) = \cup_{l=0}^\infty B_l(\Lambda)$
where $B_0(\Lambda) =\{ \emptyset \}$ the empty word.
For a word
$\mu =\mu_1\cdots \mu_k \in B_*(\Lambda)$,
 a right infinite sequence
$x =(x_i)_{i \in {\Bbb N}} \in X_\Lambda$
and
$l \in \Zp$,
put
\begin{align*}
\Gamma_l^-(\mu) & = \{ \nu_1\cdots \nu_l \in B_l(\Lambda) 
\mid \nu_1\cdots \nu_l \mu_1 \cdots \mu_k  \in B_*(\Lambda) \}, \\
\Gamma_l^-(x) & = \{ \nu_1\cdots \nu_l \in B_l(\Lambda) 
\mid (\nu_1,\cdots, \nu_l,x_1,x_2, \cdots) \in X_\Lambda \}, \\
\Gamma_l^+(\mu) & = \{ \omega_1\cdots \omega_l \in B_l(\Lambda) 
\mid  \mu_1 \cdots \mu_k \omega_1\cdots \omega_l \in B_*(\Lambda) \},\\
\Gamma_*^+(\mu) & = \cup_{l=0}^\infty \Gamma_l^+(\mu).
\end{align*}
A word
$\mu =\mu_1\cdots \mu_k \in B_*(\Lambda)$
for 
$l \in \Zp$
is said to be $l$-{\it synchronizing}\,
if 
for all $\omega \in \Gamma_*^+(\mu)$
the equality
\begin{equation*}
\Gamma_l^-(\mu) = \Gamma_l^-(\mu \omega) 
\end{equation*}
holds.
Denote by
$S_l(\Lambda)$ the set of all $l$-synchronizing words of $\Lambda$.
We say that an irreducible subshift $\Lambda$ 
is $\lambda$-{\it synchronizing}\/
if 
for any $\eta \in B_l(\Lambda)$ and $k \ge l$
there exists $\nu \in S_k(\Lambda)$ 
such that 
$\eta \nu \in S_{k-l}(\Lambda)$.
Irreducible sofic shifts are $\lambda$-synchronizing.
More generally, synchronizing subshifts are 
$\lambda$-synchronizing
(see \cite{BH} for synchronizing subshifts).
Many irreducible subshifts including Dyck shifts, $\beta$-shifts and Morse shifts 
are $\lambda$-synchronizing.
There exists a concrete  example of an irreducible subshift that is not $\lambda$-synchronizing
(see \cite{KM2010}).

\begin{prop}[\cite{2011Pre}, cf. \cite{Kr2006}, \cite{KM2010}]
$\lambda$-synchronization is invariant under not only topological
conjugacy but also flow equivalence of subshifts.
\end{prop}

For $\mu, \nu \in B_*(\Lambda)$,
we say that $\mu$ is $l$-past equivalent to $\nu$
if $\Gamma_l^-(\mu) = \Gamma_l^-(\nu)$.
We write it as  
$\mu\underset{l}{\sim}\nu$.
The following lemma is straightforward.
\begin{lem}[\cite{KM2010}, \cite{2011Pre}]
Let $\Lambda$ be a $\lambda$-synchronizing subshift.
Then we have
\begin{enumerate}
\renewcommand{\labelenumi}{(\roman{enumi})}
\item
For $\mu \in S_l(\Lambda)$,
there exists
$\mu' \in S_{l+1}(\Lambda)$
such that 
$\mu\underset{l}{\sim}\mu'$.
\item
For $\mu \in S_l(\Lambda)$,
there exist
$\beta \in \Sigma$
and
$\nu \in S_{l+1}(\Lambda)$
such that 
$\mu\underset{l}{\sim}\beta\nu$.
\end{enumerate}
\end{lem}

A $\lambda$-graph system 
is a graphical object presenting a subshift (\cite{1999DocMath}). 
It is a generalization of a finite labeled graph 
and
yields a $C^*$-algebra (\cite{2002DocMath}).  
Let ${\frak L} =(V,E,\lambda,\iota)$ be 
a $\lambda$-graph system 
 over $\Sigma$ with vertex set
$
V = \cup_{l \in \Zp} V_{l}
$
and  edge set
$
E = \cup_{l \in \Zp} E_{l,l+1}
$
with a labeling map
$\lambda: E \rightarrow \Sigma$, 
and that is supplied with  surjective maps
$
\iota( = \iota_{l,l+1}):V_{l+1} \rightarrow V_l
$
for
$
l \in  \Zp.
$
Here the vertex sets $V_{l},l \in \Zp$
are finite disjoint sets.   
Also  
$E_{l,l+1},l \in \Zp$
are finite disjoint sets.
An edge $e$ in $E_{l,l+1}$ has its source vertex $s(e)$ in $V_{l}$ 
and its terminal  vertex $t(e)$ 
in
$V_{l+1}$
respectively.
Every vertex in $V$ has a successor and  every 
vertex in $V_l$ for $l\in {\Bbb N}$ 
has a predecessor. 
It is then required that there exists an edge in $E_{l,l+1}$
with label $\alpha$ and its terminal is  $v \in V_{l+1}$
 if and only if 
 there exists an edge in $E_{l-1,l}$
with label $\alpha$ and its terminal is $\iota(v) \in V_{l}.$
For 
$u \in V_{l-1}$ and
$v \in V_{l+1},$
put
\begin{align*}
E^{\iota}_{l,l+1}(u, v)
& = \{e \in E_{l,l+1} \ | \ t(e) = v, \iota(s(e)) = u \},\\
E_{\iota}^{l-1,l}(u, v)
& = \{e \in E_{l-1,l} \ | \ s(e) = u, t(e) = \iota(v) \}.
\end{align*}
Then we require a bijective correspondence preserving their labels between 
$
E^{\iota}_{l,l+1}(u, v)
$
and
$
E_{\iota}^{l-1,l}(u, v)
$
for each pair of vertices
$u, v$.
We call this property  the local property of $\lambda$-graph system. 
We call an edge in $E$ a labeled edge and a finite sequence of connecting labeled edges a labeled path.
If a labeled path $\gamma$ labeled $\nu$
starts at a vertex $v \in V_l$ and ends at a vertex $u \in V_{l+n}$,
we say that $\nu$ leaves $v$
and write 
$s(\gamma)=v, t(\gamma) = u, \lambda(\gamma) = \nu.$  
We henceforth assume that ${\frak L}$ is left-resolving, 
which means that 
$t(e)\ne t(f)$ whenever $\lambda(e) = \lambda(f)$ for $e,f \in E$.
For a vertex
$v \in V_l$ denote by
$\Gamma_l^-(v)$ the predecessor set of $v$
which is defined by the set of  words  with length $l$
appearing as labeled paths from a vertex in $V_0$ to the vertex $v$.  
 ${\frak L}$ is said to be predecessor-separated if 
$\Gamma_l^-(v) \ne \Gamma_l^-(u)$
whenever $u, v\in V_l$  are distinct.
Two $\lambda$-graph systems
${\frak L} =(V,E,\lambda,\iota)$ over $\Sigma$
and
${\frak L}' =(V',E',\lambda',\iota')$ over $\Sigma$
are said to be isomorphic if there exist
bijections
$\varPhi_V:V \longrightarrow V'$
and 
$\varPhi_E:E \longrightarrow E'$
satisfying
$\varPhi_V(V_l) = V'_l$
and 
$\varPhi_E(E_{l,l+1}) = E'_{l,l+1}$
such that
they give rise to a labeled graph isomorphism
compatible to $\iota$ and $ \iota'$.
We note that any essential finite directed labeled graph 
${\cal G} = ({\cal V}, {\cal E}, \lambda)$ over $\Sigma$
with vertex set ${\cal V}$, 
edge set ${\cal E}$ and labeling map $\lambda:{\cal E}\longrightarrow \Sigma$
 gives rise to a $\lambda$-graph system
${\frak L}_{\cal G} =(V,E,\lambda,\iota)$
by setting
$V_l ={\cal V}, E_{l,l+1} ={\cal E}, \iota = \id$
for all $l \in \Zp$
(cf.\cite{2002DocMath}).

For a $\lambda$-synchronizing subshift $\Lambda$ over $\Sigma$,
we have introduced a $\lambda$-graph system
\begin{equation*}
{\frak L}^{\lambda(\Lambda)} 
=(V^{\lambda(\Lambda)}, E^{\lambda(\Lambda)},
\lambda^{\lambda(\Lambda)},\iota^{\lambda(\Lambda)})
\end{equation*}
defined by  $\lambda$-synchronization of $\Lambda$
as in the following way (\cite{KM2010}, \cite{2011Pre}).
Let
$V_l^{\lambda(\Lambda)}$ be the $l$-past equivalence classes of
$S_l(\Lambda)$.
We denote by $[\mu]_l$ 
the equivalence class of $\mu \in S_l(\Lambda)$.
For
$\nu \in S_{l+1}(\Lambda)$ and $\alpha \in \Gamma_1^-(\nu)$,
define a labeled edge from
$[\alpha \nu]_l\in V_l^{\lambda(\Lambda)}$ to
$[\nu]_l \in V_{l+1}^{\lambda(\Lambda)}$ labeled $\alpha$.
Such labeled edges are denoted by $E^{\lambda(\Lambda)}_{l,l+1}$.
Denote by 
$\lambda^{\lambda(\Lambda)}:E^{\lambda(\Lambda)}_{l,l+1} \longrightarrow \Sigma$
the labeling map.
Since
$S_{l+1}(\Lambda) \subset S_l(\Lambda)$,
we have a natural map
$
[\mu]_{l+1} \in V_{l+1}^{\lambda(\Lambda)} 
\longrightarrow 
[\mu]_l \in V_l^{\lambda(\Lambda)} 
$
that we denote by
$\iota^{\lambda(\Lambda)}_{l,l+1}$.
Then
${\frak L}^{\lambda(\Lambda)} 
=(V^{\lambda(\Lambda)}, E^{\lambda(\Lambda)},
\lambda^{\lambda(\Lambda)}, \iota^{\lambda(\Lambda)})
$
defines a
predecessor-separated, left-resolving 
 $\lambda$-graph stsem that presents $\Lambda$.
We call
$
{\frak L}^{\lambda(\Lambda)} 
$
the canonical
$\lambda$-synchronizing $\lambda$-graph system of $\Lambda$.

The canonical $\lambda$-synchronizing $\lambda$-graph system
may be characterized  in an intrinsic way.
Let ${\frak L} = (V, E,\lambda,\iota)$
be a predecessor-separated, left-resolving $\lambda$-graph system 
over $\Sigma$ that presents a subshift $\Lambda$.
Denote by $\{ v_1^l,\dots,v_{m(l)}^l \}$
the vertex set $V_l$ at level $l$.
For an admissible word $\nu \in B_n(\Lambda)$
and a vertex $v_i^l \in V_l$,
we say that $v_i^l$ {\it launches}\/ $\nu$
if the following two conditions hold:
\begin{enumerate}
\renewcommand{\labelenumi}{(\roman{enumi})}
\item
There exists a path labeled $\nu$ in ${\frak L}$ 
leaving the vertex $v_i^l$
and ending at a vertex in $V_{l+n}$
\item 
The word $\nu$ does not leave any other vertex in $V_l$ than $v_i^l$.
\end{enumerate}
We call the vertex $v_i^l$ the launching vertex for $\nu$.
We set
$$
S_{v_i^l}(\Lambda)
=\{ \nu \in B_*(\Lambda) \mid v_i^l \text{ launches } \nu \}.
$$

\noindent
{\bf Definition.}
A $\lambda$-graph system  ${\frak L}$
is said to be $\lambda$-{\it synchronizing} \
if for any $l \in {\Bbb N}$
and any vertex $v_i^l \in V_l$,
there exists a word $\nu \in B_*(\Lambda)$
such that 
$v_i^l$ launches $\nu$.

\begin{lem}[\cite{2011Pre}]
Keep the above notations.
Assume that ${\frak L}= (V, E,\lambda,\iota)$ 
is $\lambda$-synchronizing.
Then we have
\begin{enumerate}
\renewcommand{\labelenumi}{(\roman{enumi})}
\item
$\sqcup_{i=1}^{m(l)} S_{v_i^l}(\Lambda) = S_l(\Lambda)$.
\item
The $l$-past equivalence classes of $S_l(\Lambda)$
is $S_{v_i^l}(\Lambda), i=1,\dots,m(l)$.
\item
For any $l$-synchronizing word $w \in S_l(\Lambda)$,
there exists a vertex $v_{i(\omega)}^l \in V_l$
such that $v_{i(\omega)}^l$ launches $\omega$ 
and
$
\Gamma_l^-(\omega)= \Gamma_l^-(v_{i(\omega)}^l).
$
\end{enumerate}
\end{lem}

\noindent
{\bf Definition.}
A $\lambda$-graph system
${\frak L} = (V,E,\lambda,\iota)$ 
is said to be
$\iota$-{\it irreducible}\,
if for any two vertices $v,u \in V_l$ and 
a labeled path 
$\gamma$ starting at $u$,
there exist a labeled path 
from $v$ to a vertex $u'\in V_{l+n}$
such that
$\iota^n(u') = u$ ,
and a labeled path 
$\gamma'$ starting at $u'$
such that
$\iota^n(t(\gamma')) = t(\gamma)$
and
$\lambda(\gamma') = \lambda(\gamma)$,
where 
$t(\gamma'), t(\gamma)$
denote the terminal vertices of
$\gamma', \gamma$
respectively and
$\lambda(\gamma'),\lambda(\gamma)$
the words 
labeled  by 
$\gamma', \gamma$
respectively. 

\begin{lem}[\cite{2011Pre}]
Let
${\frak L} = (V,E,\lambda,\iota)$
be a $\lambda$-graph system
that presents a subshift $\Lambda$.
\begin{enumerate}
\renewcommand{\labelenumi}{(\roman{enumi})}
\item
If ${\frak L}$ is $\iota$-irreducible,
then
$\Lambda$ is irreducible.
\item
Assume that
${\frak L} = (V,E,\lambda,\iota)$
is $\lambda$-synchronizing.
If 
$\Lambda$ is irreducible,
then
${\frak L}$ is $\iota$-irreducible.
\end{enumerate}
\end{lem}
We then have
\begin{prop}[\cite{2011Pre}]
A subshift $\Lambda$ 
is $\lambda$-synchronizing
if and only if
there exists a left-resolving, predecessor-separated,
$\iota$-irreducible,
$\lambda$-synchronizing $\lambda$-graph system that presents $\Lambda$.
\end{prop}

\begin{thm}[\cite{2011Pre}]
For a $\lambda$-synchronizing
subshift $\Lambda$,
there uniquely exists a left-resolving, predecessor-separated,
$\iota$-irreducible,
$\lambda$-synchronizing $\lambda$-graph system that presents $\Lambda$.
The unique 
$\lambda$-synchronizing $\lambda$-graph system
is the caninical $\lambda$-synchronizing $\lambda$-graph system
${\frak L}^{\lambda(\Lambda)}$ for $\Lambda$.
\end{thm}
As in the preceding theorem,
the canonical $\lambda$-synchronizing $\lambda$-graph system
${\frak L}^{\lambda(\Lambda)}$
has a unique property in the above sense.
We  henceforth  call
${\frak L}^{\lambda(\Lambda)}$
the  $\lambda$-synchronizing $\lambda$-graph system
for $\Lambda$.

We say that a $\lambda$-graph system
${\frak L}$ is  {\it minimal}\/ if
there is no proper $\lambda$-graph subsystem of ${\frak L}$
that presents $\Lambda$.
This means that if 
${\frak L}'$ is a  $\lambda$-graph subsystem of ${\frak L}$
and 
presents the same subshift as the subshift presented by 
${\frak L}$,
then 
${\frak L}'$ coincides with ${\frak L}$.
Then we may prove that
the $\lambda$-synchronizing
$\lambda$-graph system 
${\frak L}^{\lambda(\Lambda)}$
for a $\lambda$-synchronizing subshift $\Lambda$
is minimal
(\cite{2011Pre}).

\section{$\lambda$-synchronizing $C^*$-algebras}
Let 
${\frak L} =(V,E, \lambda,\iota)$
be a left-resolving predecessor-separated $\lambda$-graph system over $\Sigma$
and 
$\Lambda$ the presented subshift by ${\frak L}$.
We denote by $\{v_1^l,\dots,v_{m(l)}^l\}$ 
the vertex set $V_l$.
Define the transition matrices $A_{l,l+1}, I_{l,l+1}$
of ${\frak L}$
by setting
for
$
i=1,2,\dots,m(l),\ j=1,2,\dots,m(l+1), \ \alpha \in \Sigma,
$ 
\begin{align*}
A_{l,l+1}(i,\alpha,j)
 & =
{\begin{cases}
1 &  
    \text{ if } \ s(e) = v_i^l, \lambda(e) = \alpha,
                       t(e) = v_j^{l+1} 
    \text{ for some }    e \in E_{l,l+1}, \\
0           & \text{ otherwise,}
\end{cases}} \\
I_{l,l+1}(i,j)
 & =
{\begin{cases}
1 &  
    \text{ if } \ \iota_{l,l+1}(v_j^{l+1}) = v_i^l, \\
0           & \text{ otherwise.}
\end{cases}} 
\end{align*}
The $C^*$-algebra $\OFL$
is realized as the universal unital $C^*$-algebra
generated by
partial isometries
$S_{\alpha}, \alpha \in \Sigma$
and projections
$E_i^l, i=1,2,\dots,m(l),\l\in \Zp 
$
 subject to the  following operator relations called $({\frak L})$:
\begin{align}
\sum_{\beta \in \Sigma} S_{\beta}S_{\beta}^*  & = 1,  \\
 \sum_{i=1}^{m(l)} E_i^l   =  1, \qquad 
 E_i^l  &  =  \sum_{j=1}^{m(l+1)}I_{l,l+1}(i,j)E_j^{l+1},  \\
 S_\alpha S_\alpha^* E_i^l & =   E_i^{l} S_\alpha S_\alpha^*, \\
S_{\alpha}^*E_i^l S_{\alpha}  =  
\sum_{j=1}^{m(l+1)}& A_{l,l+1}(i,\alpha,j)E_j^{l+1},
\end{align}
for $\alpha \in \Sigma,$
$
i=1,2,\dots,m(l),\l\in \Zp. 
$
It is nuclear (\cite[Proposition 5.6]{2002DocMath}).
For a word $\mu = \mu_1\cdots\mu_k \in B_k(X_\Lambda)$,
we set 
$S_\mu = S_{\mu_1}\cdots S_{\mu_k}.$
The algebra of all finite linear combinations of the elements of the form
$$
S_{\mu}E_i^lS_{\nu}^* \quad \text{ for }\quad \mu,\nu \in B_*(X_\Lambda), \quad i=1,\dots,m(l),\quad l \in \Zp
$$
is a dense $*$-subalgebra of  $\OFL$.
Let us denote by 
${\cal A}_{\frak L}$ the $C^*$-subalgebra of $\OFL$
generated by the projections
$E_i^l, i=1,\dots, m(l),\/ l \in \Zp$,
which is a commutative AF-algebra.
For a vertex $v_i^l \in V_l$, 
put 
\begin{align*}
\Gamma_\infty^{+}(v_i^l) 
=\{ (\alpha_1,\alpha_2,\dots, ) \in \Sigma^{\Bbb N} & \ | \
 \text{ there exists an edge } e_{n,n+1} \in E_{n,n+1} 
 \text{ for } n \ge l  \\
 \text{such that }  v_i^l =  s(e_{l,l+1}),\ &
 t(e_{n,n+1})  = s(e_{n+1,n+2}), \ 
 \lambda(e_{n,n+1}) = \alpha_{n-l+1} \}
\end{align*}
the set of all label sequences in ${\frak L}$ starting at $v_i^l$.
We say that ${\frak L}$ satisfies condition (I) if for each
$v_i^l \in V,$
the set $\Gamma_\infty^{+}(v_i^l)$ contains at least two distinct sequences.
Under the condition (I), 
the algebra $\OFL$ can be realized as the  unique 
$C^*$-algebra subject to  the relations $({\frak L})$
(\cite[Theorem 4.3]{2002DocMath}).
A $\lambda$-graph system $\frak L$ is said to $\lambda$-irreducible 
if
for an ordered pair of vertices $u, v \in V_l,$ 
there exists a number $L_l(u,v) \in {\Bbb N}$ 
such that 
for a vertex $w \in V_{l+L_l(u,v)}$ 
with
$\iota^{L_l(u,v)}(w) = u,$
there exists a path $\gamma$ in ${\frak L}$ 
such that 
$
s(\gamma) = v, \, 
t(\gamma) = w,
$
where $\iota^{L_l(u,v)}$ means the $L_l(u,v)$-times compositions of $\iota$, 
and $s(\gamma), t(\gamma)$ denote 
the source vertex, 
the terminal vertex of $\gamma$ respectively(\cite{2005MS}).
If ${\frak L}$ is $\lambda$-irreducible with  condition (I),
the $C^*$-algebra $\OFL$ is simple
(\cite[Theorem 4.7]{2002DocMath}, \cite{2005MS}).

\begin{prop}
Let
$\Lambda$ be a $\lambda$-synchronizing subshift over 
$\Sigma$ and ${\frak L}^{\lambda(\Lambda)}$
the $\lambda$-synchronizing $\lambda$-graph system for $\Lambda$.
Then $\Lambda$ is homeomorphic to a Cantor discontinuum
if and only if $\LLL$ satisfies condition (I).
\end{prop}
\begin{pf}
Assume that 
$\Lambda$ is homeomorphic to a Cantor discontinuum.
Then the right one-sided subshift $X_\Lambda$ 
is also homeomorphic to a Cantor discontinuum.
For a vertex $v_i^l \in V_l^{\lambda(\Lambda)}$,
take a $l$-synchronizing word
$\mu= \mu_1\cdots\mu_k \in S_l(\Lambda)$
such that 
$v_i^l$ launches $\mu$.
Take an infinite sequence
$x \in X_\Lambda$
such that
$\mu \in \Gamma_k^-(x)$.
Since $X_\Lambda$ 
is homeomorphic to a Cantor discontinuum,
any neighborhood of $\mu x$ in $X_\Lambda$
contains an element that is different from  $\mu x$.
Hence there exists an infinite sequence
$x' \in X_\Lambda$
such that 
$\mu x' \in X_\Lambda$
and $x \ne x'$.
As $\mu$ must leave the vertex $v_i^l$,
both the sequences 
$\mu x$ and $ \mu x'$ 
are contained in 
$\Gamma^+_\infty(v_i^l)$
so that 
$\LLL$ satisfies condition (I).

Conversely assume that 
$\LLL$ satisfies condition (I).
Since $\Lambda$ is a compact, totally disconnected metric space,
it suffices to show that  
$\Lambda$ is perfect.
It is enough to show that 
$X_\Lambda$
is perfect.
For any $x =(x_1,x_2,\dots ) \in X_\Lambda$
and a word
$\mu_1\cdots\mu_k$ 
with
$\mu_1 = x_1,\dots, \mu_k = x_k$,
consider 
a cylinder set
$U_\mu =\{ (y_n)_{n \in {\Bbb N}} \in X_\Lambda \mid y_1 = \mu_1, \dots, y_k = \mu_k \}.
$
Take an infinite path 
$(e_n)_{n \in {\Bbb N}}$ in $\LLL$ labeled $x$ 
such  that
$\lambda(e_n) = x_n, \, t(e_n) = s(e_{n+1}),\, n \in {\Bbb N}$.
Let us denote by $v_i^k \in V_k^{\lambda(\Lambda)}$
the terminal vertex of the edge 
$e_k$.
Since the follower set 
$\Gamma^+_{\infty}(v_i^k)$
of $v_i^k$ has at least two distinct sequences,
there exists
$x' = (x'_{k+1}, x'_{k+2}, \dots ) \in \Gamma^+_{\infty}(v_i^k)$
such that 
$x' \ne (x_{k+1}, x_{k+2},\dots )$.
As $x'$ starts at $v_i^k$,
the right one-sided sequence
$\mu x' =(\mu_1,\dots, \mu_k, x'_{k+1}, x'_{k+2}, \dots )$
is contained in $X_\Lambda$
and hence in $U_\mu$.
One then sees that $x$ is a cluster point in $X_\Lambda$.
\end{pf}

Let ${\frak L}= (V,E,\lambda,\iota)$
be a left-resolving, predecessor-separated $\lambda$-graph system
over $\Sigma$ that presents a $\lambda$-synchronizing subshift $\Lambda$.
Let
$S_\alpha, \alpha \in \Sigma$
and
$E_i^l, i=1,\dots,m(l), l \in \Zp$
be the generating partial isometries and the projections in
${\cal O}_{\frak L}$
satisfying the relation $({\frak L})$.
If ${\frak L}= \LLL$
the $\lambda$-synchronizing $\lambda$-graph system for $\Lambda$,
the algebra
${\cal O}_{\frak L}$
is denoted by
$\OLL$.  
We will study the algebraic structure of the $C^*$-algebra 
$\OLL$. 

\begin{lem}
If ${\frak L}$ is
the $\lambda$-synchronizing $\lambda$-graph system $\LLL$,
 we have
\begin{enumerate}
\renewcommand{\labelenumi}{(\roman{enumi})}
\item
For a vertex $v_i^l \in V_l$, 
there exists a word $\mu \in S_l(\Lambda)$ such that
$E_i^l \ge S_\mu S_\mu^*$.
\item
For a word $\mu \in S_l(\Lambda)$, 
there exists  a unique vertex $v_i^l \in V_l^{\lambda(\Lambda)}$
such that
$E_i^l \ge S_\mu S_\mu^*$.
\end{enumerate}
\end{lem}
\begin{pf}
(i)
For a vertex $v_i^l \in V_l$, 
take a word $\mu \in S_l(\Lambda)$ such that
$v_i^l$ launches $\mu$.
Since the word $\mu$ 
does not leave any other vertex in $V_l$ 
than $v_i^l$, 
we have
$S_\mu^* E_j^l S_\mu =0$ for $j \ne i$
so that
$S_\mu S_\mu^* E_j^l  =0$ for $j \ne i$.
Let $n =|\mu |$.
It then follows that
\begin{equation*}
E_i^l 
 = \sum_{\nu \in B_n(\Lambda)} S_\nu S_\nu^* E_i^l  
\ge S_\mu S_\mu^* E_i^l 
 = \sum_{j=1}^{m(l)} S_\mu S_\mu^* E_j^l 
 = S_\mu S_\mu^*.
\end{equation*}

(ii)
For a word $\mu \in S_l(\Lambda)$, 
put $v_i^l= [\mu]_l \in V_l^{\lambda(\Lambda)}$.
Since $v_i^l$ launches $\mu$,
we have 
$S_\mu^* E_j^l S_\mu =0$ for $j \ne i$
so that
$S_\mu S_\mu^* E_j^l  =0$ for $j \ne i$.
As in the above discussions,
we have
$E_i^l \ge S_\mu S_\mu^*$.
If there exists $j =1,\dots,m(l)$
such that
$E_j^l \ge S_\mu S_\mu^*$,
we have
$S_\mu^* E_j^l S_\mu  \ge S_\mu^* S_\mu \ne 0$
so that
$S_\mu^* E_j^l S_\mu \ne 0$.
Hence there exists a path in $\LLL$ labeled $\mu$
that leaves $v_j^l$.
Since $v_i^l$ launches $\mu$,
one has $j=i$.
\end{pf}
\begin{prop}
A $\lambda$-graph system  
${\frak L}$ 
is 
$\lambda$-synchronizing if and only if
for a vertex $v_i^l \in V_l$, 
there exists a word $\mu \in S_l(\Lambda)$ 
such that
$E_i^l \ge S_\mu S_\mu^*$
in $\OFL$.
\end{prop}
\begin{pf}
Since the $\lambda$-synchronizing $\lambda$-graph system for $\Lambda$
is unique and is $\LLL$,
the only if part has been proved in the preceding lemma.
We will prove the if part.
For a vertex $v_i^l \in V_l$, 
there exists a word $\mu = \mu_1\dots\mu_n \in S_l(\Lambda)$ 
such that
$E_i^l \ge S_\mu S_\mu^*$.
Hence we have
$S_\mu^* E_i^l S_\mu \ne 0$
so that 
the word $\mu $ leaves the vertex $v_i^l$
and hence $ \Gamma_l^-(v_i^l)\subset \Gamma_l^-(\mu)$.
For $\xi \in \Gamma_l^-(\mu)$
we have
$S_\xi E_i^l S_\xi^* \ge S_\xi S_\mu S_\mu^* S_\xi^* \ne 0$
so that 
$\xi \in \Gamma_l^-(v_i^l)$.
This implies
$\Gamma_l^-(\mu) \subset \Gamma_l^-(v_i^l)$
so that
\begin{equation}
 \Gamma_l^-(v_i^l) = \Gamma_l^-(\mu). \label{eqn:1} 
\end{equation}
Suppose that $\mu$ leaves $v_j^l$.
Take a path labeled $\mu$ in ${\frak L}$
from $v_j^l$ to 
$v_{j'}^{l+n} \in V_{l+n}$.
By the hypothesis,
for the vertex $v_{j'}^{l+n}$,
there exists 
$\nu \in S_{l+n}(\Lambda)$
such that
$E_{j'}^{l+n}\ge S_\nu S_\nu^*$.
By a similar argument to the above,
one knows 
\begin{equation}
 \Gamma_{l+n}^-(v_{j'}^{l+n}) = \Gamma_{l+n}^-(\nu). \label{eqn:2} 
\end{equation}
One then sees that
\begin{equation}
 \Gamma_l^-(v_j^l) = \Gamma_l^-(\mu \nu). \label{eqn:3} 
\end{equation}
One indeed sees that
for $\xi \in \Gamma_l^-(v_j^l)$,
one has
$\xi \mu \in \Gamma_{l+n}^-(v_{j'}^{l+n})$.
By (3.6), one has
$\xi \mu \in \Gamma_{l+n}^-(\nu)$
so that
$\xi \in \Gamma_l^-(\mu \nu)$.
Conversely,
for
$\eta \in \Gamma_l^-(\mu \nu)$,
one has
$\eta \mu \in \Gamma_{l+n}^-(\nu)$
so that 
by (3.6)
$\eta \mu \in \Gamma_{l+n}^-(v_{j'}^{l+n})$.
As ${\frak L}$ is left-resolving,
one has 
$\eta \in \Gamma_l^-(v_j^l)$.
Hence we have (3.7).
Now we know $\Gamma_l^-(\mu \nu) =\Gamma_l^-(\mu)$
so that  
we have
\begin{equation}
 \Gamma_l^-(v_j^l) = \Gamma_l^-(\mu ). \label{eqn:5} 
\end{equation}
By (3.5) and (3.8),
one has 
\begin{equation*}
 \Gamma_l^-(v_i^l) = \Gamma_l^-(v_j^l).  
\end{equation*}
Since ${\frak L}$ is left-resolving, one obtains that
$v_i^l = v_j^l$ and hence
$v_i^l$ launches $\mu$.
Thus  
${\frak L}$ is $\lambda$-synchronizing. 
\end{pf}
The following lemmas are stated in terms of the 
$C^*$-algebra $\OLL$ associated with the $\lambda$-synchronizing $\lambda$-graph system
$\LLL$.
\begin{lem}
For $\xi, \eta \in B_*(\Lambda)$,
we have
$\Gamma_*^+(\xi) =\Gamma_*^+(\eta)$
if and only if
$S_\xi^* S_\xi = S_\eta^* S_\eta$.
\end{lem}
\begin{pf}
Let 
$p = |\xi|, \, q = |\eta|$.
We may assume that $p \le q$.
Let $V_{t(\xi)}^p $ be the set of all terminal vertices in $V_p$ of paths in $\LLL$
labeled $\xi$, that is
$$
 V_{t(\xi)}^p = \{ v_j^p \in V_p \mid \xi \in \Gamma_p^-(v_j^p) \}.
$$
Denote by 
$\xi(p)$ the cardinal number of $V_{t(\xi)}^p$.
We write
$
V_{t(\xi)}^p = \{ v_{j_1}^p, \dots,v_{j_{\xi(p)}}^p \}.
$
Similarly,
let us denote by
$V_{t(\eta)}^q$ the set of all terminal vertices  in $V_q$ of paths in $\LLL$
labeled $\eta$.
Denote by 
$\eta(q)$ the cardinal number of $V_{t(\eta)}^q$.
We write
$
V_{t(\eta)}^q = \{ v_{k_1}^q, \dots,v_{k_{\eta(q)}}^q \}.
$
By the relation $({\frak L})$,
one sees 
\begin{equation*}
S_\xi^* S_\xi = E_{j_1}^p + \cdots + E_{j_{\xi(p)}}^p,
\qquad
S_\eta^* S_\eta = E_{k_1}^q + \cdots + E_{k_{\eta(q)}}^q.
\end{equation*}
We set
\begin{align*}
\iota^{q-p}(V_{t(\eta)}^q) & = \{ \iota^{q-p}(v_{k_1}^q), \dots,\iota^{q-p}(v_{k_{\eta(q)}}^q) \} \subset V_p,\\
\iota^{p-q}(V_{t(\xi)}^p) & = \{ v_k^q \in V_q \mid \iota^{q-p}(v_k^q) \in V_{t(\xi)}^p \}
\subset V_q.
\end{align*}
Hence we have
$S_\xi^* S_\xi = S_\eta^* S_\eta$
if and only if 
$ \iota^{p-q}(V_{t(\xi)}^p) = V_{t(\eta)}^q$.

Now assume that 
$\Gamma_*^+(\xi) =\Gamma_*^+(\eta)$.
For $v_k^q \in V_{t(\eta)}^q$,
take $\nu(k) \in S_q(\Lambda)$ such that
$v_k^q$ launches $\nu(k)$.
It is easy to see that 
$\iota^{q-p}(v_k^q)$ launches $\nu(k)$.
Since
$\nu(k) \in \Gamma_*^+(\eta)$,
one has
$\nu(k) \in \Gamma_*^+(\xi)$
so that 
$\nu(k)$ leaves a vertex in $V_{t(\xi)}^p$.
 As
$\iota^{q-p}(v_k^q)$
is the only vertex which $\nu(k)$ leaves,
one has 
$\iota^{q-p}(v_k^q) \in V_{t(\xi)}^p$.
Hence we have
$\iota^{q-p}(V_{t(\eta)}^q) \subset V_{t(\xi)}^p$
and hence
$V_{t(\eta)}^q \subset \iota^{p-q}(V_{t(\xi)}^p)$.
For the other inclusion relation,
take an arbitrary vertex 
 $v_k^p \in \iota^{p-q}(V_{t(\xi)}^p)$
and
$\mu(q) \in S_q(\Lambda)$ such that
$v_k^p$ launches $\mu(q)$.
Hence the word $\mu(q)$ leaves $\iota^{q-p}(v_k^q)$
and
$\iota^{q-p}(v_k^q)$ launches $\mu(q)$.
As $\mu(q) \in \Gamma^+_*(\xi)$,
one has 
 $\mu(q) \in \Gamma^+_*(\eta)$
so that
there exists a vertex 
$v_{k_n}^q \in V_{t(\eta)}^q$
such that 
$\mu(q)$ leaves $v_{k_n}^q$.
Therefore we have
$v_k^q = v_{k_n}^q$
and hence 
$v_k^q \in V_{t(\eta)}^q$
so that 
$\iota^{p-q}(V_{t(\xi)}^p)\subset V_{t(\eta)}^q $.
This implies
$
S_\xi^* S_\xi = S_\eta^* S_\eta.
$

Conversely assume  the equality 
$
S_\xi^* S_\xi = S_\eta^* S_\eta
$
holds
so that 
$\iota^{p-q}(V_{t(\xi)}^p) =V_{t(\eta)}^q$.
By the local property of $\lambda$-graph system,
one easily sees that 
the set of followers of $V_{t(\xi)}^p$ 
coincides with
the set of followers of $V_{t(\eta)}^q$.
This implies that
$\Gamma^+_*(\xi) =\Gamma^+_*(\eta)$. 
\end{pf}

For $\mu,\nu \in B_*(\Lambda)$,
we write
$\mu \succ \nu$
if there exists a word $\eta \in B_*(\Lambda)$ 
such that
$\Gamma_*^+(\nu) = \Gamma_*^+(\mu \eta \nu)$.
The following lemma comes from the preceding lemma.
\begin{lem}
For words $\mu,\nu \in B_*(\Lambda)$,
the following three conditions are equivalent:
\begin{enumerate}
\renewcommand{\labelenumi}{(\roman{enumi})}
\item
$\mu \succ \nu$.
\item
There exists a word $\eta \in B_*(\Lambda)$ 
such that
$ S_\nu^* S_\nu = S_\nu^* S_\eta^* S_\mu^* S_\mu S_\eta S_\nu$
in $\OLL$.
\item
There exists a word $\eta \in B_*(\Lambda)$ 
such that
$ S_\nu S_\nu^* \le S_\eta^* S_\mu^* S_\mu S_\eta$
in $\OLL$.
\end{enumerate} 
\end{lem}
\begin{pf}
The equivalence between (i) and (ii) 
follows from Lemma 3.4.
It is clear that the equality
$ S_\nu^* S_\nu = S_\nu^* S_\eta^* S_\mu^* S_\mu S_\eta S_\nu$
is equivalent to the inequality
$ S_\nu S_\nu^* \le S_\eta^* S_\mu^* S_\mu S_\eta$.
\end{pf}

\noindent
{\bf Definition.}
A $\lambda$-synchronizing subshift $\Lambda$ is said to be 
{\it synchronizingly transitive} \/ if 
for any two words $\mu,\nu \in B_*(\Lambda)$,
the both relations 
$\mu \succ \nu$ and $\nu \succ \mu$
hold.

We note that the
$\lambda$-irreduciblity
for ${\frak L}$ 
is rephrased in terms of the algebra $\OFL$ as the property that
for any $E_i^l, i=1,\dots,m(l)$,
there exists $n\in {\Bbb N}$ 
such that
$\sum_{k=1}^n \lambda_{\frak L}^k(E_i^l) \ge 1$,
where
$\lambda_{\frak L}^k(X) = \sum_{\mu \in B_k(\Lambda)}S_\mu^* X S_\mu$
for $X \in {\cal A}_{\frak L}$.

\begin{lem}
If $\Lambda$ is 
synchronizingly transitive,
then 
$\LLL$ is $\lambda$-irreducible.
\end{lem}
\begin{pf}
Take an ordered pair 
$v_i^l,v_j^l \in V_l$ of vertices.
Since $\Lambda$ is $\lambda$-synchronizing,
by Lemma 3.2
, there exists $\mu \in S_l(\Lambda)$
such that
$v_i^l$ launches $\mu$ so that
$E_i^l \ge S_\mu S_\mu^*$.
For the vertex $v_j^l$,
take a word $\nu \in B_l(\Lambda)$ 
such that $\nu \in \Gamma_l^-(v_j^l)$
so that
$S_\nu^* S_\nu \ge E_j^l$.
Now $\Lambda$ is synchronizingly transitive
so that 
we have
$$
S_\nu^* S_\eta^*S_\mu^* S_\mu S_\eta S_\nu = S_\nu^* S_\nu
$$
for some 
$\eta \in B_*(\Lambda)$,
and hence
$$
S_\nu^* S_\eta^*S_\mu^* E_i^l  S_\mu S_\eta S_\nu \ge S_\nu^* S_\nu \ge E_j^l.
$$
Put $k = | \mu \eta \nu|$. 
We then have
$
\lambda_{\LLL}^k(E_i^l) \ge E_j^l.
$
Thus we may find $n \in {\Bbb N}$ such that
$$
\sum_{k=1}^n \lambda_{\LLL}^k(E_i^l) \ge 1.
$$
\end{pf}
\begin{thm}
Let
$\Lambda$ be a $\lambda$-synchronizing subshift over 
$\Sigma$.
Assume that 
$\Lambda$ is homeomorphic to a Cantor discontinuum.
If $\Lambda$ is 
synchronizingly transitive,
then 
the $C^*$-algebra 
$\OLL$ 
associated with the 
 $\lambda$-synchronizing $\lambda$-graph system
$\LLL$ for $\Lambda$
is simple.
\end{thm}
\begin{pf}
Since 
$\Lambda$ is homeomorphic to a Cantor discontinuum,
the $\lambda$-graph system 
$\LLL$ satisfies condition (I).
By the preceding proposition,
synchronizing transitivity of $\Lambda$ implies
that $\LLL$ is $\lambda$-irreucible 
so that the $C^*$-algebra
$\OLL$ is simple by \cite[Theorem 4.7]{2002DocMath}
\end{pf}

\section{Flow equivalnce and $\lambda$-synchronizing $C^*$-algebras}

It has been proved that 
$\lambda$-synchronization is invariant under flow equivalence (\cite{2011Pre}).
In the proof, the  
Parry-Sullivan's result \cite{PS} which says that the flow equivalence relation on homeomorphisms
of Cantor sets 
is generated by topological conjugacy and expansion of symbols has been used.
Let $\Lambda$ be a subshift over alphabet 
$\Sigma =\{ 1,2,\dots, N \}$.
A new subshift $\widetilde{\Lambda}$ 
over the alphabet 
$\widetilde{\Sigma} =\{0, 1,2,\dots, N \}$
is defined 
as the subshift 
consisting of all biinfinite sequences of $\widetilde{\Sigma}$ obtained by replacing the symbol $1$ in a biinfinite sequence in the subshift $\Lambda$ by the word $01$.  
This operation is called expansion.
The Parry-Sullivan's result stated above is the following:
\begin{lem}[\cite{PS}]   
Flow equivalence relation of subshifts is generated by topological conjugacy and the expansion 
$ \Lambda \rightarrow \widetilde{\Lambda}$.
\end{lem}

In \cite{2011Pre},
it has been proved that 
the $\lambda$-synchronizing K-groups 
$K_0^\lambda (\Lambda), K_1^\lambda(\Lambda)$
 and
the $\lambda$-synchronizing Bowen-Franks groups 
$BF^0_\lambda(\Lambda),$
$ BF^1_\lambda(\Lambda)$
for 
a $\lambda$-synchronizing  subshift $\Lambda$
are invariant under flow equivalence of subshifts.
The groups 
$K_0^\lambda (\Lambda)$,
$ K_1^\lambda(\Lambda)$
 and
the Bowen-Franks groups
$BF^0_\lambda(\Lambda)$,
$ BF^1_\lambda(\Lambda)$
are realized as
the K-groups
$K_0(\OLL)$,
$K_1(\OLL)$
and
the Ext-groups
$\Ext^0(\OLL)$,
$\Ext^1(\OLL)$
for the $C^*$-algebra
$\OLL$
associated with the $\lambda$-synchronizing 
$\lambda$-graph system $\LLL$.
If the algebra
$\OLL$
is simple and purely infinite,
the K-groups
$K_0(\OLL)$,
$K_1(\OLL)$
determine the stable isomorphism class
of $\OLL$
by the structure theorem of purely infinite simple $C^*$-algebras
(\cite{Kir}, \cite{Phil}).

In this section, we will prove
that the stable isomorphism class of 
the pair
$(\OLL, \DLL)$
of $\OLL$ with its Cartan subalgebra
$\DLL$
is invariant under 
flow equivalence of $\lambda$-synchronizing subshifts.
We will not assume the simplicity of the algebra
$\OLL$.
As a result,
we also give a $C^*$-algebraic proof of the above
invariance of 
the groups 
$K_0^\lambda (\Lambda)$,
$ K_1^\lambda(\Lambda)$
 and
the Bowen-Franks groups
$BF^0_\lambda(\Lambda)$,
$ BF^1_\lambda(\Lambda)$
under flow equivalence.

Let
$\Lambda$ be a $\lambda$-synchronizing subshift over 
$\Sigma = \{ 1,2,\dots,N \}$.
Let 
$S_i,i \in \Sigma$ 
and 
$E_i^l, i=1,\dots,m(l), l \in \Zp$
be the generating partial isometries and the projections in the 
$C^*$-algebra $\OLL$ satisfying the relations $(\LLL)$.
The Cartan subalgebra
$\DLL$
is defined by the $C^*$-subalgebra of $\OLL$
generated by the projections of the form
$S_\mu E_i^l S_\mu^*, i=1,\dots,m(l), \mu \in B_*(\Lambda)$,
which is a regular maximal abelian subalgebra in $\OLL$. 
Consider the subshift
$\widetilde{\Lambda}$
 over 
$\widetilde{\Sigma}
=\{0,1,\dots,N\}$
that is obtained from $\Lambda$ 
by replacing $1$ in $\Lambda$ by $01$.
It is 
$\lambda$-synchronizing
by \cite{2011Pre}.
Denote by
$\OLTL$ the $C^*$-algebra associated with 
the 
$\lambda$-synchronizing $\lambda$-graph system
${\frak L}^{\lambda(\widetilde{\Lambda})}$
for
$\widetilde{\Lambda}$.
Similarly 
let
$\widetilde{S}_i,i \in \widetilde{\Sigma}$ 
and 
$\widetilde{E}_i^l, i=1,\dots,\tilde{m}(l), l \in \Zp$
be the generating partial isometries and the projections in the 
$C^*$-algebra $\OLTL$ satisfying the relations $(\LLTL)$.
We set the partial isometries 
\begin{equation*}
s_1 = \widetilde{S}_0\widetilde{S}_1,
\qquad
s_i = \widetilde{S}_i
\quad
\text{ for }
i=2,\dots,N
\end{equation*}
and the projection 
\begin{equation*}
P  = \widetilde{S}_0 \widetilde{S}_0^* 
   + \widetilde{S}_2 \widetilde{S}_2^* 
   + \widetilde{S}_3 \widetilde{S}_3^* + \cdots 
   + \widetilde{S}_N \widetilde{S}_N^*
   = 1 - \widetilde{S}_1 \widetilde{S}_1^*
\end{equation*}
in $\OLTL$.
\begin{lem}
$
\widetilde{S}_0^* \widetilde{S}_0 =\widetilde{S}_1 \widetilde{S}_1^* 
$
and hence
$s_1 s_1^* = \widetilde{S}_0 \widetilde{S}_0^*$,
$s_1^* s_1 = \widetilde{S}_1^* \widetilde{S}_1.$
\end{lem}
\begin{pf}
We note that 
the set $V^{\lambda(\widetilde{\Lambda})}_0$ is a singleton.
There exists a unique vertex $v_{j_0}^1 $ 
in $V^{\lambda(\widetilde{\Lambda})}_1$
such that the symbol $0$ goes to $v_{j_0}^1 $ from
$V^{\lambda(\widetilde{\Lambda})}_0$.
The vertex $v_{j_0}^1 $
is the $1$-past equivalence class 
$[1\mu]_{1}$ for a word $1 \mu \in B_*(\widetilde{\Lambda)}$.
It launches the symbol $1$.
Since $1$ is the only symbol which leaves 
$v_{j_0}^1$,
one sees
$\widetilde{S}_\alpha^* \widetilde{E}^1_{j_0}\widetilde{S}_\alpha \ne 0$
if and only if 
$\alpha = 1$.
It then follows that
\begin{equation*}
\widetilde{E}^1_{j_0} 
 = \sum_{\alpha \in \widetilde{\Sigma}}
\widetilde{S}_\alpha \widetilde{S}_\alpha^* 
\widetilde{E}^1_{j_0} 
= \widetilde{S}_1 \widetilde{S}_1^* 
\widetilde{E}^1_{j_0}. 
\end{equation*}
Hence we have
$
\widetilde{E}^1_{j_0} \le \widetilde{S}_1 \widetilde{S}_1^*. 
$
Since the  inequality
$
\widetilde{E}^1_{j_0} \ge \widetilde{S}_1 \widetilde{S}_1^* 
$
is clear,
 we have
$$
\widetilde{E}^1_{j_0} = \widetilde{S}_1 \widetilde{S}_1^*.
$$ 
As
$v_{j_0}^1 $
is the unique vertex 
in $V^{\lambda(\widetilde{\Lambda})}_1$
such that the symbol 
$0$ goes to
$v_{j_0}^1 $,
one has
$
\widetilde{S}_0^* \widetilde{S}_0 = \widetilde{E}^1_{j_0}.
$ 
The equalities
$
s_1 s_1^* 
 = \widetilde{S}_0 \widetilde{S}_0^*
$,
$
s_1^* s_1
 = \widetilde{S}_1^* \widetilde{S}_1
$
are obvious.
\end{pf}
\begin{lem}
\hspace{3cm}
\begin{enumerate}
\renewcommand{\labelenumi}{(\roman{enumi})}
\item $P = \sum_{j=1}^N s_j s_j^*$.
\item $P \ge s_\mu^* s_\mu$ for all $\mu \in B_l(\Lambda)$, $l \in {\Bbb N}$. 
\item $\sum_{\mu \in B_l(\Lambda)} s_\mu^* s_\mu \ge P$ 
for all $l \in {\Bbb N}$. 
\end{enumerate}
\end{lem}
\begin{pf}
(i)
Since $\widetilde{S}_0 \widetilde{S}_0^* = s_1 s_1^*$,
the assertion is clear.

(ii)
Since
$P = 1 - \widetilde{E}^1_{j_0}$,
it suffices to show that
$\widetilde{E}^1_{j_0}\perp s_\mu^* s_\mu$
for
$ \mu= \mu_1\cdots\mu_l \in B_l(\Lambda)$.
If $\mu_l \ne 1$,
one has  
$s_{\mu_l} =  \widetilde{S}_{\mu_l}$
so that
$
s_{\mu_l} \widetilde{S}_1 = \widetilde{S}_{\mu_l} \widetilde{S}_1 =0.
$ 
If $\mu_l = 1$,
one has  
$s_{\mu_l} =  \widetilde{S}_0 \widetilde{S}_1$
so that
$
s_{\mu_l} \widetilde{S}_1 
= \widetilde{S}_0 \widetilde{S}_1\widetilde{S}_1 =0.
$ 
In any case we have
$
s_{\mu_l} \widetilde{S}_1 = 0
$
so that
$
 s_\mu^* s_\mu \widetilde{E}^1_{j_0}  = 0.
$

(iii)
We will first prove that
$\sum_{i=1}^N s_i^* s_i \ge P$.
We know that
$\widetilde{S}_i^* \widetilde{S}_i = s_i^* s_i$ for $i=1,\dots,N$
and
$
\widetilde{S}_0^* \widetilde{S}_0 =\widetilde{S}_1 \widetilde{S}_1^* 
= 1 - P.
$
Since
$\sum_{i=0}^N \widetilde{S}_i^* \widetilde{S}_i \ge 1$
in 
$\OLTL$,
one obtains
$$
\sum_{i=0}^N \widetilde{S}_i^* \widetilde{S}_i
= 1-P + \sum_{i=1}^N s_i^* s_i \ge 1
$$
so that
$\sum_{i=1}^N s_i^* s_i \ge P$.
Suppose that
the inequality
$\sum_{\mu \in B_k(\Lambda)} s_\mu^* s_\mu \ge P$ 
holds for some $k \in {\Bbb N}$. 
It then follows that
\begin{equation*}
\sum_{\nu \in B_{k+1}(\Lambda)} s_\nu^* s_\nu
=\sum_{i=1}^N s_i^*(
\sum_{\mu \in B_k(\Lambda)} s_\mu^* s_\mu )s_i
 \ge
 \sum_{i=1}^N s_i^* P s_i
= \sum_{i,j=1}^N s_i^*  s_j  s_j^* s_i
= \sum_{i=1}^N s_i^*  s_i
\ge P.
\end{equation*}
Hence we have the desired inequalities.
\end{pf}
In the $\lambda$-graph system $\LLL$,
recall that the set 
$\Gamma_l^-(v_i^l)$ 
for  a vertex $v_i^l$ in $V_l$
denotes the predecessor of
$v_i^l$ which is the set of words
of $B_l(\Lambda)$ presenting by labeled paths 
terminating at $v_i^l$.  
Put the projections
for $i=1,2,\dots,m(l), \, l \in \Zp$
\begin{equation*}
e_i^l = 
\prod_{\mu \in \Gamma_l^-(v_i^l)}s_\mu^* s_\mu 
\prod_{\nu \in B_l(\Lambda)\backslash\Gamma_l^-(v_i^l)}(P-s_\nu^* s_\nu). 
\end{equation*}
For 
$\mu \in B_*(\Lambda)$,
put
$$
s_\mu^* s_\mu^{1} =  s_\mu^* s_\mu,
\qquad
s_\mu^* s_\mu^{-1} =  P - s_\mu^* s_\mu.
$$
For 
$v_i^l \in V_l^{\lambda(\Lambda)}$,
define a function
$f_i^l: B_l(\Lambda) \longrightarrow \{1,-1\}$
by setting
\begin{equation*}
f_i^l(\mu) =
\begin{cases}
 1 & \text{ if } \mu \in  \Gamma_l^-(v_i^l),\\
-1 & \text{ if } \mu \not\in  \Gamma_l^-(v_i^l)
\end{cases}
\end{equation*}
so that
$$
e_i^l = 
\prod_{\mu \in B_l(\Lambda)}s_\mu^* s_\mu^{f_i^l(\mu)}. 
$$
Denote by 
$\{1, -1\}^{B_l(\Lambda)}$
the set of all functions from 
$B_l(\Lambda)$ to $\{ 1, -1\}$. 
\begin{lem}
For
$
\epsilon \in 
\{1, -1\}^{B_l(\Lambda)}$, 
we have
$
\prod_{\mu \in B_l(\Lambda)}s_\mu^* s_\mu^{\epsilon(\mu)} \ne 0
$
if and only if
$\epsilon = f_i^l$ for some
$i=1,\dots,m(l)$. 
In this case
$
\prod_{\mu \in B_l(\Lambda)}s_\mu^* s_\mu^{\epsilon(\mu)} = e_i^l.
$
\end{lem}
\begin{pf}
Suppose that
$\epsilon = f_i^l$ for some
$i=1,\dots,m(l)$. 
Since $\Lambda$ 
is $\lambda$-synchronizing,
there exists $\nu \in S_l(\Lambda)$ such that
$v_i^l$ launches $\nu$
so that
\begin{align*}
s_\mu^* s_\mu 
& \ge s_\nu s_\nu^* \quad \text{ for } \quad \mu \in \Gamma_l^-(v_i^l),\\
P -s_\mu^* s_\mu 
& \ge s_\nu s_\nu^* \quad \text{ for } \quad \mu \in 
B_l(\Lambda)\backslash\Gamma_l^-(v_i^l).
\end{align*}
Hence we have
$
\prod_{\mu \in B_l(\Lambda)}s_\mu^* s_\mu^{f_i^l(\mu)} \ge s_\nu s_\nu^* \ne 0.$

Conversely suppose that
$
\prod_{\mu \in B_l(\Lambda)}s_\mu^* s_\mu^{\epsilon(\mu)} \ne 0.
$
Since
$
\prod_{\mu \in B_l(\Lambda)}s_\mu^* s_\mu^{\epsilon(\mu)} 
\in {\cal A}_{\lambda(\widetilde{\Lambda})},
$
there exists
$k \ge l$ and $i_1 = 1,2,\dots,\tilde{m}(k)$
such that
$
\prod_{\mu \in B_l(\Lambda)}s_\mu^* s_\mu^{\epsilon(\mu)} 
\ge
\widetilde{E}_{i_1}^k  
\in {\cal A}_{\lambda(\widetilde{\Lambda})}.
$
Take 
$\omega \in S_k(\widetilde{\Lambda})$
such that
$v_{i_1}^k$ launches $\omega$.
Since
$\sum_{\mu \in B_l(\Lambda)}s_\mu^*s_\mu \ge P$,
there exists 
$\mu \in B_l(\Lambda)$ such that
$s_\mu^* s_\mu 
\ge 
\widetilde{E}_{i_1}^k$.
Hence we see that
$\mu \omega \in B_*(\Lambda)$.
As the rightmost of $\mu$ is not $0$,
the leftmost of $\omega$  is not $1$.
Let
$\bar{\omega}$
be
the word 
in
$B_*(\Lambda)$
obtained from 
$\omega$ by putting $1$ in place of $01$ in $\omega$.
Since
$\widetilde{E}_{i_1}^k \ge \widetilde{S}_\omega \widetilde{S}_\omega^*$,
one sees that
$$
\prod_{\mu \in B_l(\Lambda)}s_\mu^* s_\mu^{\epsilon(\mu)} 
\ge
s_{\bar{\omega}}s_{\bar{\omega}}^*.
$$
As
$[ \bar{\omega}]_l \in V_l^{\lambda(\Lambda)}$,
we have
 $[ \bar{\omega}]_l = v_i^l$
 for some $i=1,\dots,m(l)$.
 The vertex $v_i^l$ launches $\bar{\omega}$
 so that 
 $\epsilon = f_i^l$.
\end{pf}

\begin{lem}
For $\mu, \nu \in B_l(\Lambda)$ and $\alpha, \beta \in \Sigma$,
we have
\begin{enumerate}
\renewcommand{\labelenumi}{(\roman{enumi})}
\item $s_\mu^*(P - s_\alpha^* s_\alpha )s_\mu\cdot s_\mu^* s_\beta^* s_\beta s_\mu
= (P - s_{\alpha\mu}^* s_{\alpha\mu} ) s_{\beta\mu}^* s_{\beta\mu}$.
\item $s_\alpha^* \cdot s_\mu^*s_\mu (P - s_\nu^*s_\nu)s_\alpha
= s_{\mu\alpha}^*  s_{\mu\alpha} (P - s_{\nu\alpha}^*s_{\nu\alpha})$. 
\end{enumerate}
\end{lem}
\begin{pf}
(i)
Since $P  s_\beta^* s_\beta =  s_\beta^* s_\beta$
and hence
$
s_\mu^*P  s_\beta^* s_\beta s_\mu
 =  s_{\beta\mu}^* s_{\beta\mu}$,
we have
\begin{align*}
s_\mu^*(P - s_\alpha^* s_\alpha )s_\mu\cdot s_\mu^* s_\beta^* s_\beta s_\mu
&=s_\mu^* P s_\beta^* s_\beta s_\mu - s_\mu^* s_\alpha^* s_\alpha s_\beta^* s_\beta s_\mu \\
&=P s_{\beta\mu}^* s_{\beta\mu} - s_{\alpha\mu}^* s_{\alpha\mu} s_{\beta\mu}^* s_{\beta\mu}\\
&=(P - s_{\alpha\mu}^* s_{\alpha\mu} ) s_{\beta\mu}^* s_{\beta\mu}.  
\end{align*}

(ii)
Since 
$P s_\alpha  = s_\alpha$
and
$
s_{\mu \alpha}^* s_{\mu \alpha} = s_{\mu \alpha}^* s_{\mu \alpha} P$,
we have
\begin{align*}
s_\alpha^* \cdot s_\mu^*s_\mu (P - s_\nu^*s_\nu)s_\alpha
& = s_{\mu\alpha}^*  s_{\mu\alpha} - s_{\mu\alpha}^* s_{\mu\alpha}s_{\nu\alpha}^*s_{\nu\alpha} \\
& = s_{\mu\alpha}^*  s_{\mu\alpha} (P - s_{\nu\alpha}^*s_{\nu\alpha}).
\end{align*}
\end{pf}

\begin{lem}
The partial isometries
$s_{\alpha}, \alpha \in \Sigma$
and 
the projections
$e_i^l, i=1,2,\dots,m(l)$,
$l\in \Zp 
$
satisfy the  following operator relations:
\begin{align}
\sum_{\beta \in \Sigma} s_{\beta}s_{\beta}^*  & = P,  \\
 \sum_{i=1}^{m(l)} e_i^l   =  P, \qquad 
 e_i^l  &  =  \sum_{j=1}^{m(l+1)}I_{l,l+1}(i,j)e_j^{l+1},  \\
 s_\alpha s_\alpha^*e_i^l & =   e_i^{l} s_\alpha s_\alpha^*, \\
s_{\alpha}^*e_i^l s_{\alpha}  =  
\sum_{j=1}^{m(l+1)}& A_{l,l+1}(i,\alpha,j)e_j^{l+1},
\end{align}
for $\alpha \in \Sigma,$
$
i=1,2,\dots,m(l),\l\in \Zp, 
$
where
$I_{l,l+1}, A_{l,l+1}$ denote the transition mattices
for the $\lambda$-graph system $\LLL$.
\end{lem}
\begin{pf}
The equality (4.1) has been proved in Lemma 4.3 (i).

It follows that
\begin{equation*}
P =
\prod_{\mu \in B_l(\Lambda)}(s_\mu^*s_\mu + P - s_\mu^*s_\mu)
=\sum_{\epsilon \in \{ -1,1\}^{B_l(\Lambda)}}
\prod_{\mu \in B_l(\Lambda)}s_\mu^*s_\mu^{\epsilon(\mu)}.
\end{equation*} 
By Lemma 4.4, 
the nonzero 
$\prod_{\mu \in B_l(\Lambda)} s_\mu^*s_\mu^{\epsilon(\mu)}$
is of the form
$\prod_{\mu \in B_l(\Lambda)} s_\mu^*s_\mu^{f_i^l(\mu)}$
for some
$i=1,\dots,m(l)$
so that we have
$P = \sum_{i=1}^{m(l)}e_i^l.$

We will next show the equality (4.4).
It follows that
\begin{align*}
s_\alpha^* e_i^l s_\alpha
& = s_\alpha^*(
\prod_{\mu \in \Gamma_l^-(v_i^l)}s_\mu^* s_\mu 
\prod_{\nu \in B_l(\Lambda)\backslash\Gamma_l^-(v_i^l)}(P-s_\nu^* s_\nu) 
) s_\alpha \\
& = 
\prod_{\mu \in \Gamma_l^-(v_i^l)}s_{\mu\alpha}^* s_{\mu\alpha} 
\prod_{\nu \in B_l(\Lambda)\backslash\Gamma_l^-(v_i^l)}(P-s_{\nu\alpha}^* s_{\nu\alpha}). 
\end{align*}
Hence 
$s_\alpha^* e_i^l s_\alpha$
is written as a finite sum of 
$e_j^{l+1}, j=1,\dots,m(l+1)$.
If
$s_\alpha^* e_i^l s_\alpha
\ge e_j^{l+1}$,
then one has
\begin{align*}
s_\alpha^*(s_\mu^* s_\mu)  s_\alpha
& \ge e_j^{l+1} \quad \text{ for } \mu \in \Gamma_l^-(v_i^l),  \\
s_\alpha^*(P - s_\nu^* s_\nu)  s_\alpha
& \ge e_j^{l+1} \quad \text{ for } \nu \in B_l(\Lambda)\backslash\Gamma_l^-(v_i^l) . 
\end{align*}
Since
\begin{equation*}
e_j^{l+1}
 =
\prod_{\xi \in \Gamma_{l+1}^-(v_j^{l+1})}s_\xi^* s_\xi 
\prod_{\eta\in B_{l+1}(\Lambda)\backslash\Gamma_{l+1}^-(v_j^{l+1})}
(P-s_\eta^* s_\eta) 
\end{equation*}
and 
$\Lambda$ is $\lambda$-synchronizing,
there exists 
$\zeta(j) \in S_{l+1}(\Lambda)$
such that
$[\zeta(j)]_{l+1} = v_j^{l+1}$.
Hence we have
$e_j^{l+1} \ge s_{\zeta(j)}s_{\zeta(j)}^*$.
As
$s_\alpha^* e_i^l s_\alpha \ge e_j^{l+1}
\ge s_{\zeta(j)}s_{\zeta(j)}^*$,
one has
$e_i^l \ge s_{\alpha\zeta(j)}s_{\alpha\zeta(j)}^* \ne 0$.
Hence
\begin{align*}
\mu\alpha\zeta(j) \in B_*(\Lambda)
& \quad \text{ for } \mu \in \Gamma_l^-(v_i^l), \\
\nu\alpha\zeta(j) \not\in B_*(\Lambda)
&  \quad \text{ for } \nu \in B_l(\Lambda)\backslash\Gamma_l^-(v_i^l) 
\end{align*}
so that
$[\alpha \zeta(j)]_l = v_i^l$.
Since
$[\zeta(j)]_{l+1} = v_j^{l+1}$,
one has
$A_{l,l+1}(i,\alpha,j) =1$.
Therefore 
the condition
$s_\alpha^* e_i^l s_\alpha \ge e_j^{l+1}$
implies
$A_{l,l+1}(i,\alpha,j) =1$.
Hence we obtain
$$
s_\alpha^* e_i^l s_\alpha  
= \sum_{j=1}^{m(l+1)}A_{l,l+1}(i,\alpha,j)e_j^{l+1}.
$$

We will next prove the second equality of (4.2).
By the equalities 
\begin{align*}
 e_i^l 
& = 
\prod_{\mu \in \Gamma_l^-(v_i^l)}s_\mu^* s_\mu 
\prod_{\nu \in B_l(\Lambda)\backslash\Gamma_l^-(v_i^l)}(P-s_\nu^* s_\nu) 
 \\
& = 
\prod_{\mu \in \Gamma_l^-(v_i^l)}
(\sum_{k=1}^{m(1)}s_\mu^* e_k^1 s_\mu ) 
\prod_{\nu \in B_l(\Lambda)\backslash\Gamma_l^-(v_i^l)}
(P - \sum_{h=1}^{m(1)}s_\nu^*e_h^1 s_\nu) 
\end{align*}
one knows that
$e_i^l$ is a finite sum of $e_1^l,\dots, e_{m(l+1)}^{l+1}$.
Suppose that
 $e_i^l \ge e_j^{l+1}$.
Since $v_j^{l+1} = [\zeta(j)]_{l+1}$
for some $\zeta(j) \in S_{l+1}(\Lambda)$,
one has
$e_j^{l+1}\ge s_{\zeta(j)}s_{\zeta(j)}^*$
and hence
$e_i^l\ge s_{\zeta(j)}s_{\zeta(j)}^*$.
This implies
$$
\prod_{\mu \in \Gamma_l^-(v_i^l)}s_\mu^* s_\mu 
\prod_{\nu \in B_l(\Lambda)\backslash\Gamma_l^-(v_i^l)}(P-s_\nu^* s_\nu) 
\ge 
s_{\zeta(j)}s_{\zeta(j)}^*
$$
so that
\begin{align*}
s_\mu^* s_\mu \ge s_{\zeta(j)}s_{\zeta(j)}^*
\text{ and hence } s_{\mu\zeta(j)} \ne 0
& \quad \text{ for } \mu \in \Gamma_l^-(v_i^l), \\
P - s_\nu^* s_\nu \ge s_{\zeta(j)}s_{\zeta(j)}^*
\text{ and hence } s_{\nu\zeta(j)} = 0
&  \quad \text{ for } \nu \in B_l(\Lambda)\backslash\Gamma_l^-(v_i^l). 
\end{align*}
Hence
\begin{align*}
\mu \zeta(j) \in B_*(\Lambda) 
& \quad \text{ for } \mu \in \Gamma_l^-(v_i^l), \\
\nu\zeta(j) \not\in B_*(\Lambda)
&  \quad \text{ for } \nu \in B_l(\Lambda)\backslash\Gamma_l^-(v_i^l). 
\end{align*}
Thus 
one has
$[\zeta(j)]_l = v_i^l$.
As
$[\zeta(j)]_{l+1} = v_j^{l+1}$,
one obtaines that
$I_{l,l+1}(i,j) =1$.
We then conclude 
the second equality of (4.2).

The projections
$e_i^l$ and $s_\alpha^*s_\alpha$
all belong to the 
commutative $C^*$-subalgebra 
 of 
$\OLTL$
generated by the projections
$
\widetilde{S}_\mu \widetilde{S}_{\xi_1}^*\widetilde{S}_{\xi_1}\cdots 
\widetilde{S}_{\xi_n}^*\widetilde{S}_{\xi_n}
\widetilde{S}_{\mu}^*,
\mu, \xi_1\cdots\xi_n \in B_*(\widetilde{\Lambda}).
$
The commutativity between 
$e_i^l$ and $s_\alpha^*s_\alpha$
is obvious. 
Thus we complete the proof.
\end{pf}

Therefore we have
\begin{cor}
The $C^*$-subalgebra of $\OLTL$
generated by the partial isometries 
$s_\alpha, \alpha \in \Sigma$
and the projections
$e_i^l, i=1,\dots,m(l), \, l\in \Zp$
is canonically isomorphic to the 
$C^*$-algebra 
$\OLL$ associated to the $\lambda$-graph system $\LLL$.
\end{cor}
We identify the algebra 
$\OLL$ with the above
$C^*$-subalgebra of $\OLTL$
generated by the partial isometries 
$s_\alpha, \alpha \in \Sigma$
and the projections
$e_i^l, i=1,\dots,m(l), \, l\in \Zp$.
We note that
the projections 
$e_i^l, i=1,\dots,m(l), \, l\in \Zp$
and
$P$ are written by  
$s_\alpha, s_\alpha^*, \alpha \in \Sigma$
so that
the subalgebra
$\OLL$
is generated by 
$s_\alpha, \alpha \in \Sigma$.

We will henceforce prove that 
the $C^*$-subalgebra
$P\OLTL P$ is generated by 
$s_\alpha, \alpha \in \Sigma$,
that is
$P\OLTL P = \OLL$.
Let
$\ALTL$
be the $C^*$-subalgebra
of
$\OLTL$ generated by the projections
$\widetilde{E}_i^l$,
$ i=1,\dots,\tilde{m}(l), l \in \Zp$,
similarly 
$\ALL$
the $C^*$-subalgebra
of
$\OLTL$ generated by the projections
$e_i^l$,
$ i=1,\dots,\tilde{m}(l), l \in \Zp$.
The subalgebra 
$\ALL$ is naturally regarded as a corresponding subalgebra of
$\OLL$ through the canonical isomorphism in the above corollary.

For a word 
$\nu=\nu_1\cdots\nu_l  
\in B_l(\widetilde{\Lambda})$
satisfying
$\nu_1 \ne 1, \nu_l \ne 0$,
we define
the word 
$\bar{\nu}
\in
B_*(\Lambda)$ 
by
putting $1$ in place of $01$ in $\nu$.
Since $s_1 =\widetilde{S}_0 \widetilde{S}_1$,
the following lemma is straightforward. 
\begin{lem}
For any $\mu = \mu_1 \cdots\mu_k \in B_k(\widetilde{\Lambda})$,
the partial isometry
$\widetilde{S}_\mu$ is of the form:
\begin{equation*}
\widetilde{S}_\mu
=
\begin{cases}
 s_{\bar{\mu}}  
 & \text{ if } \mu_1 \ne 1, \ \mu_k \ne 0,\\
 \widetilde{S}_1 s_{\overline{\mu_2\cdots\mu_k}}  
 & \text{ if } \mu_1 = 1, \ \mu_k \ne 0,\\
 s_{\overline{\mu_1\cdots\mu_{k-1}}} \widetilde{S}_0 
 & \text{ if } \mu_1 \ne 1, \ \mu_k = 0,\\
 \widetilde{S}_1 s_{\overline{\mu_2\cdots\mu_{k-1}}} \widetilde{S}_0 
 & \text{ if } \mu_1 = 1, \ \mu_k = 0.
\end{cases}
\end{equation*}
\end{lem}
\begin{lem}
For any $\mu = \mu_1 \cdots\mu_k \in B_k(\widetilde{\Lambda})$,
we have
\begin{equation*}
\widetilde{S}_\mu P
=
\begin{cases}
 s_{\bar{\mu}} P & \text{ if } \mu_1 \ne 1, \ \mu_k \ne 0,\\
 \widetilde{S}_1 s_{\overline{\mu_2\cdots\mu_k}} P 
 & \text{ if } \mu_1 = 1, \ \mu_k \ne 0,\\
 0 
 & \text{ if } \mu_1 \ne 1, \ \mu_k = 0,\\
 0 
 & \text{ if } \mu_1 = 1, \ \mu_k = 0.
\end{cases}
\end{equation*}
\end{lem}
\begin{pf}
By the preceding lemma,
it suffices to show that
$\widetilde{S}_0 P =0$ for both the third case and the fourth case.
As 
$\widetilde{S}_0^*\widetilde{S}_0 = \widetilde{S}_1\widetilde{S}_1^*$,
we have
\begin{equation*}
\widetilde{S}_0 P 
= \widetilde{S}_0 \widetilde{S}_1\widetilde{S}_1^*P
= \widetilde{S}_0 \widetilde{S}_1\widetilde{S}_1^*
(1 - \widetilde{S}_1\widetilde{S}_1^*)
=0.
\end{equation*}
\end{pf}

\begin{lem}
For any $\mu = \mu_1 \cdots\mu_k \in B_k(\widetilde{\Lambda})$,
we have
\begin{equation*}
P\widetilde{S}_\mu^*\widetilde{S}_\mu P
=
\begin{cases}
P s_{\bar{\mu}}^* s_{\bar{\mu}} P 
& \text{ if } \mu_1 \ne 1, \ \mu_k \ne 0,\\
P s_{\overline{\mu_2\cdots\mu_k}}^*s_1^* s_1 s_{\overline{\mu_2\cdots\mu_k}} P 
& \text{ if } \mu_1 = 1, \ \mu_k \ne 0,\\
0 
& \text{ if } \mu_1 \ne 1, \ \mu_k = 0,\\
0 
& \text{ if } \mu_1 = 1, \ \mu_k = 0.
\end{cases}
\end{equation*}
\end{lem}
\begin{pf}
By the preceding lemma,
it suffices to show the equality for the second case.
For $\mu_1 =1, \mu_k \ne 0$,
one has
$\widetilde{S}_\mu P = \widetilde{S}_1 s_{\overline{\mu_2\cdots\mu_k}} P$
so that 
$$
P\widetilde{S}_\mu^*\widetilde{S}_\mu P
=P s_{\overline{\mu_2\cdots\mu_k}}^* \widetilde{S}_1^* \widetilde{S}_1 s_{\overline{\mu_2\cdots\mu_k}} P 
=P s_{\overline{\mu_2\cdots\mu_k}}^* s_1^* s_1 s_{\overline{\mu_2\cdots\mu_k}} P. 
$$
\end{pf}
\begin{cor}
$P \ALTL P = \ALL.$
\end{cor}
\begin{pf}
By the previous lemma,
one sees that
for $\mu \in B_*(\widetilde{\Lambda})$,
the element
$P\widetilde{S}_\mu^*\widetilde{S}_\mu P$
belongs to $P\ALL P$.
As $P$ is the unit of $\ALL$,
one knows that 
$P\widetilde{S}_\mu^*\widetilde{S}_\mu P \in \ALL$.
Since
$\ALTL$ is generated by the projections
$\widetilde{S}_\mu^*\widetilde{S}_\mu, \mu \in B_*(\widetilde{\Lambda})$,
we have
$P \ALTL P \subset \ALL$.
The converse inclusion relation
$P \ALTL P \supset \ALL$
is clealr.
\end{pf}
\begin{lem}
For any $\mu = \mu_1 \cdots\mu_k \in B_k(\widetilde{\Lambda})$,
we have
\begin{equation*}
(1-P)\widetilde{S}_\mu^*\widetilde{S}_\mu (1-P)
=
\begin{cases}
 \widetilde{S}_1 s_{\overline{\mu_1\cdots\mu_k 1}}^* 
s_{\overline{\mu_1\cdots\mu_k 1}} \widetilde{S}_1^*
& \text{ if } \mu_1 \ne 1,\\
\widetilde{S}_1 s_{\overline{\mu_2\cdots\mu_k 1}}^* s_1^* 
s_1 s_{\overline{\mu_2\cdots\mu_k 1}} \widetilde{S}_1^*
& \text{ if } \mu_1 = 1.
\end{cases}
\end{equation*}
\end{lem}
\begin{pf}
Since
$1 -P = \widetilde{S}_1 \widetilde{S}_1 ^*$,
it follows that
\begin{equation*}
(1-P)\widetilde{S}_\mu^*\widetilde{S}_\mu (1-P)
=\widetilde{S}_1 \widetilde{S}_{\mu 1}^*\widetilde{S}_{\mu 1}\widetilde{S}_1 ^*
= 
\begin{cases}
 \widetilde{S}_1 s_{\overline{\mu_1\cdots\mu_k 1}}^* 
s_{\overline{\mu_1\cdots\mu_k 1}} \widetilde{S}_1^*
& \text{ if } \mu_1 \ne 1,\\
\widetilde{S}_1 s_{\overline{\mu_2\cdots\mu_k 1}}^* \widetilde{S}_1^* 
\widetilde{S}_1 s_{\overline{\mu_2\cdots\mu_k 1}} \widetilde{S}_1^*
& \text{ if } \mu_1 = 1.
\end{cases}
\end{equation*}
As
$\widetilde{S}_1^*\widetilde{S}_1 = s_1^* s_1$,
one sees the desired equalities.
\end{pf}
\begin{cor}
$(1-P) \ALTL (1- P) \subset \widetilde{S}_1 \ALL \widetilde{S}_1^*.$
\end{cor}
\begin{pf}
By the previous lemma,
one sees that
for $\mu \in B_*(\widetilde{\Lambda})$,
the element
$(1-P)\widetilde{S}_\mu^*\widetilde{S}_\mu (1-P)$
belongs to $\widetilde{S}_1 \ALL \widetilde{S}_1^*$
so that we have
$(1-P) \ALTL (1- P) \subset \widetilde{S}_1 \ALL \widetilde{S}_1^*.$
\end{pf}
\begin{prop}
$P \OLTL P \subset \OLL.$
\end{prop}
\begin{pf}
The $C^*$-algebra 
$P \OLTL P$
is generated by the elements of the form:
$$
P \widetilde{S}_\mu \widetilde{S}_{\xi_1}^* \widetilde{S}_{\xi_1}\cdots 
  \widetilde{S}_{\xi_n}^* \widetilde{S}_{\xi_n}
  \widetilde{S}_\nu^* P,
  \qquad
  \mu,\xi_1,\dots,\xi_n, \nu \in B_*(\widetilde{\Lambda}).
$$
Suppose that
$ 
P \widetilde{S}_\mu \widetilde{S}_{\xi_1}^* \widetilde{S}_{\xi_1}\cdots 
  \widetilde{S}_{\xi_n}^* \widetilde{S}_{\xi_n}
  \widetilde{S}_\nu^* P \ne 0.
$
Let
$\mu = \mu_1\cdots\mu_k, \nu = \nu_1\cdots\nu_h$.
Since
$
P \widetilde{S}_\mu = \widetilde{S}_\mu \ne 0
$
and
$ 
\widetilde{S}_\nu^* P =\widetilde{S}_\nu^*\ne 0,
$
one has
$\mu_1 \ne 1, \nu_1\ne 1$.
Hence the words 
$\mu, \nu$ 
satisfy the first condition or the third condition in Lemma 4.8.

Case 1: $\mu_k \ne 0, \nu_h \ne 0$.

Since
$ 
\widetilde{S}_{\mu_k} \widetilde{S}_1\widetilde{S}_1^* =0,
$
we have
$ 
\widetilde{S}_{\mu_k}(1 - P)  =0
$
so that 
$ 
\widetilde{S}_\mu P =\widetilde{S}_\mu.
$
Hence 
$ 
\widetilde{S}_\mu
$
commutes with $P$.
Similarly 
$ 
\widetilde{S}_\nu
$
commutes with $P$.
By Lemma 4.8,
one sees that
$ 
\widetilde{S}_\mu = s_{\overline{\mu}},
\widetilde{S}_\nu = s_{\overline{\nu}}.
$
It then follows that
\begin{equation*}
P \widetilde{S}_\mu \widetilde{S}_{\xi_1}^* \widetilde{S}_{\xi_1}\cdots 
  \widetilde{S}_{\xi_n}^* \widetilde{S}_{\xi_n}
  \widetilde{S}_\nu^* P
 = s_{\overline{\mu}} P  \widetilde{S}_{\xi_1}^* \widetilde{S}_{\xi_1}\cdots 
  \widetilde{S}_{\xi_n}^* \widetilde{S}_{\xi_n} P s_{\overline{\nu}}^*.
\end{equation*}
Since
$\widetilde{S}_{\xi_1}^* \widetilde{S}_{\xi_1}\cdots 
  \widetilde{S}_{\xi_n}^* \widetilde{S}_{\xi_n} \in \ALTL
$
and
$P\ALTL P = \ALL$,
one sees that the element
\begin{equation*}
P \widetilde{S}_\mu \widetilde{S}_{\xi_1}^* \widetilde{S}_{\xi_1}\cdots 
  \widetilde{S}_{\xi_n}^* \widetilde{S}_{\xi_n}
  \widetilde{S}_\nu^* P
\end{equation*}
belongs to
$
  s_{\overline{\mu}} \ALL s_{\overline{\nu}}^*
$
and hence to
$\OLL$.

Case 2: $\mu_k \ne 0, \nu_h = 0$.

As in the above discussion,
we know that 
$ 
\widetilde{S}_\mu
$
commutes with $P$.
Since 
$P \widetilde{S}_0^*\widetilde{S}_0 =0$,
one has
\begin{align*}
P \widetilde{S}_\mu \widetilde{S}_{\xi_1}^* \widetilde{S}_{\xi_1}\cdots 
  \widetilde{S}_{\xi_n}^* \widetilde{S}_{\xi_n}
  \widetilde{S}_\nu^* P
& = 
  \widetilde{S}_\mu P \widetilde{S}_{\xi_1}^* \widetilde{S}_{\xi_1}\cdots 
  \widetilde{S}_{\xi_n}^* \widetilde{S}_{\xi_n}
  \widetilde{S}_0^* \widetilde{S}_{\nu_1\cdots \nu_{h-1}}^* P \\
& = 
  \widetilde{S}_\mu P \widetilde{S}_0^*\widetilde{S}_0
  \widetilde{S}_{\xi_1}^* \widetilde{S}_{\xi_1}\cdots 
  \widetilde{S}_{\xi_n}^* \widetilde{S}_{\xi_n}
  \widetilde{S}_0^* \widetilde{S}_{\nu_1\cdots \nu_{h-1}}^* P = 0
\end{align*}
a contradiction.

Case 3: $\mu_k = 0, \nu_h \ne 0$.

This case is similar to Case 2.

Case 4: $\mu_k = 0, \nu_h = 0$.

Since
$\widetilde{S}_0 P =0$,
we have
$\widetilde{S}_\mu = \widetilde{S}_\mu (1 - P)$
and similarly
$\widetilde{S}_\nu^* = (1 -P)\widetilde{S}_\nu^*$.
As both words $\mu, \nu$ satisfy the third condition in Lemma 4.8,
one sees that
$$
\widetilde{S}_\mu = s_{\overline{\mu_1\cdots\mu_{k-1}}}\widetilde{S}_0,
\qquad
\widetilde{S}_\nu = s_{\overline{\nu_1\cdots\nu_{h-1}}}\widetilde{S}_0.
$$
It then follows that
$$
P\widetilde{S}_\mu =\widetilde{S}_\mu 
= s_{\overline{\mu_1\cdots\mu_{k-1}}}\widetilde{S}_0(1 - P),
\qquad
\widetilde{S}_\nu^* P  = \widetilde{S}_\nu^* 
= (1 - P) \widetilde{S}_0^*s_{\overline{\nu_1\cdots\nu_{h-1}}}^*.
$$
Hence we have
\begin{align*}
& P \widetilde{S}_\mu \widetilde{S}_{\xi_1}^* \widetilde{S}_{\xi_1}\cdots 
  \widetilde{S}_{\xi_n}^* \widetilde{S}_{\xi_n}
  \widetilde{S}_\nu^* P \\
= & s_{\overline{\mu_1\cdots\mu_{k-1}}}\widetilde{S}_0(1 - P)
\widetilde{S}_{\xi_1}^* \widetilde{S}_{\xi_1}\cdots 
  \widetilde{S}_{\xi_n}^* \widetilde{S}_{\xi_n}
(1 - P) \widetilde{S}_0^*s_{\overline{\nu_1\cdots\nu_{h-1}}}^*.
\end{align*}
By the preceding lemma,
one knows that
$(1-P) \ALTL (1 - P) \subset \widetilde{S}_1 \ALL \widetilde{S}_1^* 
$
so that the element 
$
\widetilde{S}_0(1 - P)
\widetilde{S}_{\xi_1}^* \widetilde{S}_{\xi_1}\cdots 
  \widetilde{S}_{\xi_n}^* \widetilde{S}_{\xi_n}
(1 - P) \widetilde{S}_0^*
$
belongs to 
$
\widetilde{S}_0
\widetilde{S}_1 \ALL \widetilde{S}_1^*
\widetilde{S}_0^*
$
which is 
$s_1 \ALL s_1^*$.
Hence  
the element
$
P \widetilde{S}_\mu \widetilde{S}_{\xi_1}^* \widetilde{S}_{\xi_1}\cdots 
  \widetilde{S}_{\xi_n}^* \widetilde{S}_{\xi_n}
  \widetilde{S}_\nu^* P
$
belongs to 
$s_1 \ALL s_1^*$
and hence to 
$\OLL$.

Therefore in all cases we see that 
$
P \widetilde{S}_\mu \widetilde{S}_{\xi_1}^* \widetilde{S}_{\xi_1}\cdots 
  \widetilde{S}_{\xi_n}^* \widetilde{S}_{\xi_n}
  \widetilde{S}_\nu^* P
$
belongs to
$\OLL$ so that
we conclude
$
P\OLTL P \subset \OLL$.
\end{pf}
Let
${\cal D}_{\lambda(\widetilde{\Lambda})}$
be the $C^*$-subalgebra
of
$\OLTL$ generated by the projections
$\widetilde{S}_\mu\widetilde{E}_i^l\widetilde{S}_\mu^*$,
$\mu\in B_*(\widetilde{\Lambda})$,
$ i=1,\dots,\tilde{m}(l), l \in \Zp$,
similarly 
${\cal D}_{\lambda(\Lambda)}$
the $C^*$-subalgebra
of
$\OLTL$ generated by the projections
$s_\nu e_i^l s_\nu^*$,
$\nu\in B_*(\Lambda)$,
$i=1,\dots,\tilde{m}(l), l \in \Zp$.
The subalgebra 
${\cal D}_{\lambda(\Lambda)}$ 
is naturally regarded as a corresponding subalgebra of
$\OLL$ through the canonical isomorphism in Corollary 4.7.
\begin{prop}
\hspace{3cm}
\begin{enumerate}
\renewcommand{\labelenumi}{(\roman{enumi})}
\item $P \OLTL P = \OLL.$
\item $ \OLTL P \OLTL = \OLTL.$
\item $ P \DLTL P = \DLL.$
\end{enumerate}
\end{prop}
\begin{pf}
(i)
The inclusion relation
$P \OLTL P\supset  \OLL$
is obvious so that by the preceding proposition
we have
$P \OLTL P = \OLL.$

(ii)
Since
$
\widetilde{S}_0^*
\widetilde{S}_0 = 
\widetilde{S}_1
\widetilde{S}_1^*
$
one has
$
\widetilde{S}_0^* P \widetilde{S}_0 
=
\widetilde{S}_0^*
\widetilde{S}_0 
= 
\widetilde{S}_1
\widetilde{S}_1^*.
$
It follows that
$$
\widetilde{S}_0^* P \widetilde{S}_0  + P
=\sum_{j=0}^N
\widetilde{S}_j \widetilde{S}_j^* =1.
$$
This means that $P$ is a full projection in $\OLTL$.

(iii)
In the proof of Proposition 4.14,
the projection 
$
P \widetilde{S}_\mu \widetilde{S}_{\xi_1}^* \widetilde{S}_{\xi_1}\cdots 
  \widetilde{S}_{\xi_n}^* \widetilde{S}_{\xi_n}
  \widetilde{S}_\mu^* P
$
belongs to
$\DLL$ so that
$
P\DLTL P \subset \DLL$.
The other inclusion relation
$
P\DLTL P \supset \DLL$
is clear.
\end{pf}
Let
$K(H)$ be the $C^*$-algebra of all compact operators on a separable infinite dimensional Hilbert space $H$
and
$C(H)$ a maximal commutative $C^*$-subalgebra of $K(H)$.
\begin{thm}
Assume that $\Lambda$ is a $\lambda$-synchronizing subshift
that is homeomorphic to a Cantor discontinuum.
Then we have
$$
(\OLTL \otimes K(H), \DLTL \otimes C(H) ) \cong (\OLL \otimes K(H), \DLL \otimes C(H)).
$$ 
In particular we have
$$
\OLTL \otimes K(H) \cong \OLL \otimes K(H).
$$ 
\end{thm}
\begin{pf}
Proposition 4.15 (ii) 
shows that the projection $P$ is full in $\OLTL$.
By \cite{Bro}, we have desired assertions.  
\end{pf}
Therefore we conclude 
\begin{thm}
Assume that $\lambda$-synchronizing subshifts 
$\Lambda_1$ and $\Lambda_2$
are homeomorphic to a Cantor discontinuum.
Supposee that 
$\Lambda_1$ is flow equivalent to $\Lambda_2$.
Then we have
$$
({\cal O}_{\lambda({\Lambda_1})}  \otimes K(H), 
{\cal D}_{\lambda({\Lambda_1})}   \otimes C(H) )
\cong 
({\cal O}_{\lambda({\Lambda_2})}  \otimes K(H), 
{\cal D}_{\lambda({\Lambda_2})}   \otimes C(H)). 
$$ 
In particular we have
$$
{\cal O}_{\lambda({\Lambda_1})}  \otimes K(H) 
\cong
{\cal O}_{\lambda({\Lambda_2})}  \otimes K(H).
$$ 
\end{thm}
\begin{pf}
Flow equivalence relation of subshifts is generated by topological conjugacy and expansion
$\Lambda \longrightarrow \widetilde{\Lambda}$.
Suppose that $\lambda$-synchronizing subshifts 
$\Lambda_1$ and $\Lambda_2$ are topologically conjugate.
By \cite[Proposition 3.5]{KM2010},
their symbolic matrix systems 
$(\M^{\lambda(\Lambda_1)}, I^{\lambda(\Lambda_1)})
$
and 
$(\M^{\lambda(\Lambda_2)}, I^{\lambda(\Lambda_2)})
$
are strong shift equivalence.
Then we have
$$
({\cal O}_{\lambda({\Lambda_1})}  \otimes K(H), 
{\cal D}_{\lambda({\Lambda_1})}   \otimes C(H) )
\cong 
({\cal O}_{\lambda({\Lambda_2})}  \otimes K(H), 
{\cal D}_{\lambda({\Lambda_2})}   \otimes C(H)). 
$$ 
by \cite[Theorem 4.4]{2004ETDS}.
Hence by the above theorem, we have desired assertions.
\end{pf}

\begin{cor}[\cite{2011Pre}]
Assume that $\lambda$-synchronizing subshifts 
$\Lambda_1$ and $\Lambda_2$
are homeomorphic to a Cantor discontinuum.
Supposee that 
$\Lambda_1$ is flow equivalent to $\Lambda_2$.
Then 
the $\lambda$-synchronizing K-groups
and
the $\lambda$-synchronizing Bowen-Franks groups 
are isomorphic to each other, that is 
$$
K_i^\lambda (\Lambda_1)\cong  K_i^\lambda(\Lambda_2)
\quad
\text{ and }
\quad
BF^i_\lambda(\Lambda_1)\cong BF^i_\lambda(\Lambda_2),
\qquad i=0,1.
$$
\end{cor}

\section{Examples}
{\bf 1. Sofic shifts}.

Let $\Lambda$ be an irreducible sofic shift which is homeomorphic to a Cantor discontinuum.
Let $\G_{F(\Lambda)}$
be a finite directed labeled graph of the minimal left-resolving presentation of $\Lambda$.
Such a labeled graph is unique up to graph isomorphism and called the left Fischer cover
(\cite{Fis}, \cite{Kr84}, \cite{Kr87}, \cite{We}).
Let
${\frak L}_{\G_{F(\Lambda)}}$
be the $\lambda$-graph system
associated to the finite labeled graph $\G_{F(\Lambda)}$
(see \cite[Proposition 8.2]{2002DocMath}).
Then the $\lambda$-synchronizing $\lambda$-graph system
$\LLL$ for the sofic shift $\Lambda$
is nothing but the $\lambda$-graph system
${\frak L}_{\G_{F(\Lambda)}}$.
Let $N$ 
be the number of the vertices of the graph $\G_{F(\Lambda)}$.
Let $\M_{F(\Lambda)}$
be the 
$N \times N$ symbolic matrix of the graph $\G_{F(\Lambda)}$.
Let  $A_{F(\Lambda)}$
be the $N \times N$ nonnegative matrix 
defined from $\M_{F(\Lambda)}$ by all symbols equal to $1$
in each component of $\M_{F(\Lambda)}$.
Then 
the $C^*$-algebra 
$\OLL$ of the $\lambda$-graph system $\LLL$ is simple, purely infinite. 
The algebra 
$\OLL$ is nothing but the labeled graph $C^*$-algebra
${\cal O}_{\G_{F(\Lambda)}}$
for the labeled graph $\G_{F(\Lambda)}$ (cf. \cite{BP}).
It is isomorphic to a Cuntz-Krieger algebra.
The $\lambda$-synchronizing $K$-groups and the Bowen-Franks groups 
are as follows:
\begin{align*}
K_0^\lambda (\Lambda) & = {\Bbb Z}^N/(I_N - A_{F(\Lambda)}^t){\Bbb Z}^N,
\qquad 
 K_1^\lambda(\Lambda) = \Ker(I_N - A_{F(\Lambda)}^t) \quad \text{ in } {\Bbb Z}^N \\
\intertext{ and }
BF^0_\lambda(\Lambda)& = {\Bbb Z}^N/(I_N - A_{F(\Lambda)}){\Bbb Z}^N,
\qquad 
BF^1_\lambda(\Lambda) = \Ker(I_N - A_{F(\Lambda)}) \quad \text{ in } {\Bbb Z}^N.
\end{align*}
They are all invariant under flow equivalence of $\Lambda$.

{\bf 2. Dyck shifts.}

Let $N > 1$ be a fixed positive integer. 
We consider the Dyck shift $D_N$ 
with alphabet 
$\Sigma = \Sigma^- \cup \Sigma^+$
where
$\Sigma^- = \{ \alpha_1,\dots,\alpha_N \},
\Sigma^+ = \{ \beta_1,\dots,\beta_N \}.
$
The symbols 
$\alpha_i, \beta_i$
correspond to 
the brackets
$(_i,  )_i$
respectively.
The Dyck inverse monoid for $\Sigma$  has the relations
\begin{equation}
\alpha_i \beta_j
=
\begin{cases}
 {\bold 1} & \text{ if } i=j,\\
 0  & \text{ otherwise} 
\end{cases}  
\end{equation}
for 
$ i,j = 1,\dots,N$ (\cite{Kr}, \cite{KM2003}).
A word 
$\omega_1\cdots\omega_n$ 
of $\Sigma$
is admissible for $D_N$ 
precisely if
$
\prod_{m=1}^{n} \omega_m \ne 0.
$
For a word $\omega= \omega_1 \cdots \omega_n $ of $\Sigma,$ 
we denote by $\tilde{\omega}$ 
its reduced form.
Namely $\tilde{\omega}$ is a word of 
$\Sigma \cup \{ 0, {\bold 1} \}$
obtained after the operations (5.1).
Hence a word $\omega$ of $\Sigma$
is forbidden for $D_N$ if 
and only if $\tilde{\omega} = 0$.

Let us describe the Cantor horizon $\lambda$-graph system ${\frak L}^{Ch(D_N)}$ of $D_N$.
Let $\Sigma_N$ be the full $N$-shift 
$\{ 1,\dots,N \}^{\Bbb Z}$.
We denote by 
$B_l(D_N)$
and 
$B_l(\Sigma_N)$
 the set of admissible words of length 
$l$ of $D_N$
and that of 
$\Sigma_N$ respectively.
The vertices $V_l$ of ${\frak L}^{Ch(D_N)}$
 at level $l$
are given by the words of length $l$
consisting of the symbols of $\Sigma^+$.
That is, 
$$
V_l = \{ \beta_{\mu_1}\cdots\beta_{\mu_l} \in B_l(D_N)
 \mid \mu_1\cdots\mu_l\in B_l(\Sigma_N) \}.
$$
It is easy to see that each word of $V_l$ is $l$-synchronizing in $D_N$
such that
$V_l$ represent the all $l$-past equivalence classes of $D_N$.
Hence
we know that 
$V_l = V^{\lambda(D_N)}_l.$
The cardinal number of $V_l$ is $N^l$.
The mapping $\iota ( = \iota_{l,l+1}) :V_{l+1}\rightarrow V_l$ 
deletes the rightmost symbol of a word such as 
\begin{equation} 
\iota(\beta_{\mu_1}\cdots\beta_{\mu_{l+1}}) 
= \beta_{\mu_1}\cdots\beta_{\mu_l},
\qquad \beta_{\mu_1}\cdots\beta_{\mu_{l+1}} \in V_{l+1}.
\end{equation}
There exists an edge labeled $\alpha_j$ 
from
$\beta_{\mu_1}\cdots\beta_{\mu_l} \in V_l$ 
to
$\beta_{\mu_0}\beta_{\mu_1}\cdots\beta_{\mu_l} \in V_{l+1}$
precisely if
$\mu_0 = j,$
and 
there exists an edge labeled $\beta_j$ 
from
$\beta_j\beta_{\mu_1}\cdots\beta_{\mu_{l-1}} \in V_l$ 
to
$\beta_{\mu_1}\cdots\beta_{\mu_{l+1}} \in V_{l+1}.$
The resulting labeled Bratteli diagram with 
$\iota$-map
is the Cantor horizon $\lambda$-graph system ${\frak L}^{Ch(D_N)}$ of $D_N$.
One knows easily the following:
\begin{prop}
The Dyck shift $D_N$ is $\lambda$-synchronizing, 
and the $\lambda$-synchronizing $\lambda$-graph system
${\frak L}^{\lambda(D_N)}$
is
the Cantor horizon $\lambda$-graph system ${\frak L}^{Ch(D_N)}$.
\end{prop}
The Cantor horizon $\lambda$-graph system 
${\frak L}^{Ch(D_N)}$
 gives rise to a purely infinite simple $C^*$-algebra
${\cal O}_{{\frak L}^{Ch(D_N)}}$ (\cite{KM2003},\cite{2007JOT}).
The K-groups of the  $C^*$-algebra
${\cal O}_{{\frak L}^{Ch(D_N)}}$
are realized as the K-groups of the $\lambda$-graph system 
${\frak L}^{Ch(D_N)}$ 
so that  
$$
K_0({\cal O}_{\lambda(D_N)}) 
\cong {\Bbb Z}/N{\Bbb Z} \oplus C({\frak K},{\Bbb Z}), \qquad 
K_1({\cal O}_{\lambda(D_N}) \cong 0 \qquad (\cite{KM2003},\cite{2007JOT}),
$$
where 
$
C({\frak K},{\Bbb Z})
$
denotes the abelian group of all ${\Bbb Z}$-valued continuous functions 
on a Cantor discontinuum ${\frak K}$.

{\bf 3. Topological Markov Dyck shifts.}

We  consider a generalization of 
the above discussions for the Dyck shifts.
Let $A =[A(i,j)]_{i,j=1,\dots,N}$
be an $N\times N$ matrix with entries in $\{0,1\}$.
Consider
the Dyck inverse monoid 
for the alphabet
$\Sigma = \Sigma^- \cup \Sigma^+$
where
$\Sigma^- = \{ \alpha_1,\cdots,\alpha_N \},
\Sigma^+ = \{ \beta_1,\cdots,\beta_N \}.
$
which  has the relations (5.1).
Let ${\cal O}_A$ 
be
the Cuntz-Krieger algebra for the matrix $A$
that is the universal $C^*$-algebra generated by 
$N$ partial isometries $t_1,\dots,t_N$ 
subject to the following relations:
$$
\sum_{j=1}^N t_j t_j^* = 1, 
\qquad
t_i^* t_i = \sum_{j=1}^N A(i,j) t_jt_j^* \quad \text{ for } i= 1,\dots,N
$$
(\cite{CK}).
Define a correspondence 
$\varphi_A :\Sigma \longrightarrow \{t_i^*, t_i \mid i=1,\dots,N\}$
by setting
$$ 
\varphi_A(\alpha_i) = t_i^*,\qquad 
\varphi_A(\beta_i) = t_i,  \qquad i=1,\dots,N.
$$
We denote by $\Sigma^*$ the set of all words 
$\gamma_1\cdots \gamma_n$ of elements of $\Sigma$.
Define the set
$$
{\frak F}_A = \{ \gamma_1\cdots \gamma_n \in \Sigma^* \mid
\varphi_A(\gamma_1)\cdots \varphi_A( \gamma_n) = 0 \text{ in } {\cal O}_A \}.
$$
Let $D_A$ be the subshift over $\Sigma$ whose forbidden words are 
${\frak F}_A.$
The subshift is called the topological Markov Dyck shift defined by $A$
(\cite{2010MS}).
These kinds of  subshifts have first appeared in \cite{Kr2006BLM} 
in semigroup setting
and in \cite{HIK} 
in more general setting without using $C^*$-algebras (cf. \cite{2010MS}).
If all entries of $A$  are $1$, the subshift becomes the Dyck shift $D_N$
with $2N$ bracket, because  
the partial isometries 
$\{ \varphi_A(\alpha_i), \varphi(\beta_i) \mid i=1,\dots,N\}$
yield the  Dyck inverse monoid.
Consider the following  subsystem of $D_A$
\begin{equation*}
D_A^+  = \{ {(\gamma_i)}_{i \in \Bbb Z} \in D_A \mid
\gamma_i \in \Sigma^+, i \in \Bbb Z \}, 
\end{equation*}
which is identified with the topological Markov shift 
$$
\Lambda_A = \{ {(x_i)}_{i \in \Bbb Z}\in \{ 1,\dots,N \}^{\Bbb Z}
 \mid A(x_i,x_{i+1}) = 1, i \in \Bbb Z \}
$$ 
defined by the matrix $A$.
If $A$ satisfies condition (I) in the sense of Cuntz-Krieger \cite{CK},
the subshift $D_A$ is not sofic (\cite[Proposition 2.1]{2010MS}.
Hence most irreducible matrix $A$ yield  non Markov subshifts $D_A$.
Similarly to the Dyck shifts,
one may consider
the Cantor horizon $\lambda$-graph systems 
$\LCHDA$
for the topological Markov Dyck  shifts $D_A$,
which
have been studied in \cite{2010MS}.
We denote by 
$B_l(D_A^+)$
 the set of admissible words of length 
$l$ of $D_A^+$.
The vertices $V_l$, $l \in \Zp$
 of  $\LCHDA$ 
are given by the admissible words of length $l$
consisting of the symbols of $\Sigma^+$.
They are $l$-synchronizing words of $D_A$
such that
the $l$-past equivalence classes of them
coincide with the $l$-past equivalence classes of 
the set of all $l$-synchronizing words of $D_A$.
Hence $V_l = V^{\lambda(D_A)}_l$.
Since $V_l$ is identified with $B_l(\Lambda_A)$,
we may write $V_l$ as
$$
V_l = \{ \beta_{\mu_1} \cdots \beta_{\mu_l}
 \mid \mu_1\cdots\mu_l\in B_l(\Lambda_A) \}.
$$
The mapping $\iota ( = \iota_{l,l+1}) :V_{l+1}\rightarrow V_l$ 
is defined by deleting the rightmost symbol of a corresponding word as in (5.2).
There exists an edge labeled $\alpha_j$ from
$\beta_{\mu_1}\cdots \beta_{\mu_l}\in V_l$ 
to
$\beta_{\mu_1}\cdots \beta_{\mu_{l+1}} \in V_{l+1}$
precisely if
$\mu_0 = j,$
and 
there exists an edge labeled $\beta_j$ 
from
$\beta_j\beta_{\mu_1}\cdots \beta_{\mu_{l-1}}
\in V_l$ 
to
$\beta_{\mu_1}\cdots \beta_{\mu_{l+1}}$.
It is easy to see that the resulting labeled Bratteli diagram with 
$\iota$-map
becomes a $\lambda$-graph system
written 
$\LCHDA$
called 
the Cantor horizon $\lambda$-graph system 
for the topological Markov Dyck  shifts $D_A$. 
\begin{prop}
The subshift $D_A$ is $\lambda$-synchronizing,
and  the $\lambda$-synchronizing $\lambda$-graph system
${\frak L}^{\lambda(D_A)}$
is
the Cantor horizon $\lambda$-graph system
${\frak L}^{Ch(D_A)}$.
\end{prop}
Hence the $C^*$-algebra 
${\cal O}_{\lambda(D_A)}$ coincides with 
is the algebra
${\cal O}_{{\frak L}^{Ch(D_A)}}$.
By \cite[Lemma 2.5]{2010MS},
if $A$ satisfies condition (I) in the sense of \cite{CK},
the $\lambda$-graph system $\LCHLA$ satisfies $\lambda$-condition (I)
in the sense of \cite{2005MS}.
 If $A$ is irreducible,
the $\lambda$-graph system $\LCHLA$ is $\lambda$-irreducible.
Hence we have

\begin{prop}
Suppose that  $A$ is an irreducible matrix with entries in $\{0,1\}$ satisfying condition (I). 
Then the $C^*$-algebra 
${\cal O}_{\lambda(D_A)}$ associated with 
the $\lambda$-synchronizing  $\lambda$-graph system
${\frak L}^{\lambda(D_A)}$ for the topological Markov Dyck shift $D_A$
is simple and  purely infinite.
\end{prop}

One knows that 
$\beta$-shifts for $1 <\beta\in {\Bbb R}$,
a schynchronizing counter shift named as the context free shift in \cite[Example 1.2.9]{LM},
and Motzkin shifts 
are all $\lambda$-synchronizing.
Their $C^*$-algebras for the $\lambda$-synchronizing $\lambda$-graph systems
have been studied in the papers
 \cite{KMW},  \cite{1999JOT}, \cite{2004MZ}
respectively.

\medskip

{\it Acknowledgment:}
The author would like to thank Wolfgang Krieger
for his various discussions  and constant encouragements.


\end{document}